\documentclass[letterpaper]{article}
\usepackage{uai2020}
\usepackage[margin=1in]{geometry}

\usepackage{tikz}
\usetikzlibrary{arrows}
\usetikzlibrary{positioning}
\usetikzlibrary{arrows.meta}
\usepackage{enumitem}

\usepackage{amsthm}
\usepackage{amsmath}
\usepackage{amssymb}  
\usepackage[round]{natbib}

\usepackage{caption}
\usepackage{subcaption}
\usepackage{graphicx}
\usepackage{xr}
\usepackage{color}
\definecolor{bluegray}{rgb}{0.04,0,0.7}
\definecolor{darkbrown}{rgb}{0.40,0.2,0.05}
\usepackage[colorlinks,citecolor=bluegray,linkcolor=darkbrown,urlcolor=blue,breaklinks]{hyperref}

\usepackage{tabularx}
\usepackage{color}
\usepackage{setspace}

\usepackage[ruled,linesnumbered]{algorithm2e}
\SetKwInOut{Input}{input}
\SetKwInOut{Output}{output}
\SetKw{Throw}{throw}

\newtheorem{theorem}{Theorem}[section]
\newtheorem{corollary}[theorem]{Corollary}
\newtheorem{lemma}[theorem]{Lemma}
\newtheorem{definition}[theorem]{Definition}
\newtheorem{example}[theorem]{Example}
\newtheorem{proposition}[theorem]{Proposition}

\theoremstyle{definition}

 \newenvironment{proofof}[1][]{\begin{trivlist}
 \item[\hskip \labelsep {\bfseries Proof of #1.}]}{\hfill{}$\square$\end{trivlist}}

\theoremstyle{remark}

\definecolor{blue-violet}{rgb}{0.54, 0.17, 0.89}
\definecolor{antiquefuchsia}{rgb}{0.57, 0.36, 0.51}
\definecolor{amethyst}{rgb}{0.6, 0.4, 0.8}
\definecolor{blue-violet}{rgb}{0.54, 0.17, 0.89}
\definecolor{ao}{rgb}{0.0, 0.5, 0.0}
\definecolor{blue(ncs)}{rgb}{0.0, 0.53, 0.74}
\definecolor{dgreen}{rgb}{0.12, 0.3, 0.17}
\definecolor{cadmiumgreen}{rgb}{0.0, 0.42, 0.24}
\definecolor{darkolivegreen}{rgb}{0.33, 0.42, 0.18}
\definecolor{dartmouthgreen}{rgb}{0.05, 0.5, 0.06}

\newcommand{\tild}{\raise.17ex\hbox{ $\scriptstyle\sim$ }}

\DeclareMathOperator{\Var}{Var}

\DeclareMathOperator{\Cov}{Cov}
\DeclareMathOperator{\PDAG}{PDAG}
\DeclareMathOperator{\DAG}{DAG}
\DeclareMathOperator{\MPDAG}{MPDAG}
\DeclareMathOperator{\PossDe}{PossDe}

\DeclareMathOperator{\De}{De}

\DeclareMathOperator{\An}{An}
\DeclareMathOperator{\PossAn}{PossAn}
\DeclareMathOperator{\Pa}{Pa}
\DeclareMathOperator{\pa}{pa}

\DeclareMathOperator{\Forbb}{Forb}

\newcommand{\mb}[1]{\mathbf{#1}}

\newcommand{\mpdag}{$\MPDAG$}
\newcommand{\Mpdag}{$\MPDAG$}

\newcommand{\vars}[1][V]{\mathbf{#1}}
\newcommand{\e}[1][E]{\mathbf{#1}}

\newcommand{\g}[1][G]{\mathcal{#1}}

\newcommand{\f}[2][X,Y]{\Forbb(#1,#2)}
\newcommand{\fb}[2][X,Y]{\Forbb(\mathbf{#1},#2)}

\newcommand{\pstar}[1][p]{{#1}^{*}}

\usepackage{times}

\title{Identifying causal effects in  \\ maximally oriented partially directed acyclic graphs}

\author{ {\bf Emilija Perkovi\'c} \\
Department of Statistics \\
University of Washington\\
Seattle, WA 98195-4322
}

\begin{document}
\maketitle

\begin{abstract}
We develop a necessary and sufficient causal identification criterion for maximally oriented partially directed acyclic graphs (MPDAGs). MPDAGs as a class of graphs include directed acyclic graphs (DAGs), completed partially directed acyclic graphs (CPDAGs), and CPDAGs with added background knowledge. As such, they represent the type of graph that can be learned from observational data and background knowledge under the assumption of no latent variables. 
Our identification criterion can be seen as a generalization of the g-formula of \cite{robins1986new}.
We further obtain a generalization of  the truncated factorization formula  \citep{Pearl2009} and compare our criterion to the generalized adjustment criterion of \cite{perkovic17} which is sufficient, but not necessary for causal identification.
\end{abstract}

\section{INTRODUCTION} \label{sec:intro}

The gold standard method for answering causal questions are randomized controlled trials. 
 In some cases, however, it may be impossible, unethical, or simply too expensive to perform a desired experiment. For this purpose, it is of interest to consider whether a causal effect can be identified from observational data. 
 
We consider the problem of  identifying causal effects from a causal graph that represents the observational data under the assumption of causal sufficiency. 
If the causal directed acyclic graph (DAGs, \citealp[e.g.][]{Pearl2009}) is known, then all causal effects can be identified and estimated from observational data (see e.g.\ \citealp{robins1986new,pearl1995causal,pearl1995probabilistic,galles1995testing}).

In general, however, it is not possible to learn the underlying causal DAG from observational data.
 When all variables in the causal system are observed, one can at most learn a completed partially directed acyclic graph (CPDAG,  \citealp{meek1995causal,andersson1997characterization,spirtes2000causation,Chickering02}). A CPDAG uniquely represents a Markov equivalence class of DAGs (see Section \ref{sec:prelims} for definitions). 

If in addition to observational data one has background knowledge of some pairwise causal relationships, one can obtain a maximally oriented partially directed acyclic graph (\mpdag{})  which uniquely represents a refinement of the Markov equivalence class of DAGs \citep{meek1995causal}. Other types of background knowledge, such as tiered orderings, data from previous experiments, or specific model restrictions also induce \mpdag{}s \citep{tetrad1998,hoyer08,hauserBuehlmann12,eigenmann17,wang2017permutation,ernestroth2016}. 

\begin{figure}
\tikzstyle{every edge}=[draw,>=stealth',->]
\newcommand\dagvariant[1]{\begin{tikzpicture}[xscale=.5,yscale=0.5]
\node (a) at (0,0) {};
\node (d) at (0,2) {};
\node (b) at (2,0) {};
\node (c) at (2,2) {};
\begin{scope}[gray]
\draw (a) edge [-] (b);
\draw (b) edge [-] (c);
\draw (d) edge [-] (a);
\draw (d) edge [-] (b);
\draw (a) edge [-] (c);
\end{scope}
\draw #1;
\end{tikzpicture}}

\centering
\begin{subfigure}{.4\columnwidth}
  \centering
\begin{tikzpicture}[->,>=latex,shorten >=1pt,auto,node distance=.8cm,scale=1,transform shape]
  \tikzstyle{state}=[inner sep=1pt, minimum size=12pt]

  \node[state] (Xia) at (0,0) {$Y_1$};
  \node[state] (Xka) at (0,2) {$V_1$};
  \node[state] (Xja) at (2,0) {$X$};
  \node[state] (Xsa) at (2,2) {$Y_2$};

 \draw 		(Xka) edge [line width=1.1pt,-] (Xja);
 \draw 		(Xia) edge  [-] (Xka);
 \draw    	(Xsa)  edge [-] (Xja);
 \draw    	(Xja) edge [-] (Xia);
 \draw    	(Xsa) edge [-] (Xia);
\end{tikzpicture}
  \caption{}
  \label{cpdag11}
\end{subfigure}
\unskip
\vrule
\hspace{.5cm}
\begin{subfigure}{.50\columnwidth}
\dagvariant{
(a)  edge [->]  (b)
(b)  edge [->] (c)
(d) edge [->] (a)
(d) edge [->] (b)
(a) edge [->] (c)
}
\dagvariant{
(a)  edge [->]  (b)
(b)  edge [->]  (c)
(a) edge [->] (d)
(d) edge [->] (b)
(a) edge [->] (c)
}
\dagvariant{
(a)  edge [->]  (b)
(b)  edge [->]  (c)
(a) edge [->] (d)
(b) edge [->] (d)
(a) edge [->] (c)
}
\dagvariant{
(b) edge [->] (a)
(b) edge [->] (c)
(a) edge [->] (d)
(b) edge [->] (d)
(a) edge [->] (c)
}
\dagvariant{
(b) edge [->] (a)
(b) edge [->] (c)
(d) edge [->] (a)
(b) edge [->] (d)
(a) edge [->] (c)
}
\dagvariant{
(b) edge [->] (a)
(b) edge [->] (c)
(d) edge [->] (a)
(d) edge [->] (b)
(a) edge [->] (c)
}
\dagvariant{
(a) edge [->] (b)
(c) edge [->] (b)
(a) edge [->] (d)
(b) edge [->] (d)
(a) edge [->] (c)
}
\dagvariant{
(b) edge [->] (a)
(c) edge [->] (b)
(a) edge [->] (d)
(b) edge [->] (d)
(c) edge [->] (a)
}
\dagvariant{
(b) edge [->] (a)
(b) edge [->] (c)
(a) edge [->] (d)
(b) edge [->] (d)
(c) edge [->] (a)
}
\dagvariant{
(a) edge [->] (b)
(c) edge [->] (b)
(a) edge [->] (d)
(b) edge [->] (d)
(c) edge [->] (a)
}
\dagvariant{
(a) edge [->] (b)
(b) edge [->] (c)
(a) edge [->] (d)
(d) edge [->] (b)
(c) edge [->] (a)
}
\caption{}
\label{alldags-cpdag}
\end{subfigure}
\caption{(a) CPDAG $\g[C]$, (b) DAGs represented by $\g[C]$.}
\label{fig-cpdag-dag-fig}
\end{figure}

To understand the difference and connections between DAGs, CPDAGs and MPDAGs, consider graphs in Figures \ref{fig-cpdag-dag-fig} and \ref{fig:example_meekgraph11}. 
Graph $\g[C]$ in Figure \ref{fig-cpdag-dag-fig}(a) is an example of a CPDAG that can be learned given enough observational data on variables $X, V_1, Y_1,$ and $ Y_2$. All DAGs in the Markov equivalence class represented by $\g[C]$ are given in Figure \ref{fig-cpdag-dag-fig}(b).
Graph $\g$ in Figure \ref{fig:example_meekgraph11}(a) is an MPDAG that can be obtained from CPDAG $\g[C]$ in Figure  \ref{fig-cpdag-dag-fig}(a) and background knowledge that $Y_1$ is a cause of $X$ and that $X$ is a cause of $Y_2$ (see \citealp{meek1995causal} for details on incorporating this type of background knowledge). 
All DAGs represented by $\g$ are given in Figure \ref{fig:example_meekgraph11}(b) and are a subset of DAGs in Figure \ref{fig-cpdag-dag-fig}(b).

One can consider \mpdag{}s as a graph class that is generally more causally  informative than CPDAGs and  less causally informative than DAGs. Conversely,  a CPDAG can be seen as a special case of an \mpdag{} when the added background knowledge is not additionally informative compared to the observational data. Similarly, a DAG is a special case of an \mpdag{} when the additional background knowledge is fully causally informative. We will use  \mpdag{}s to refer to all graphs in this paper. 

The topic of identifying causal effects in  \mpdag{}s  has generated a wealth of research in recent years. The most relevant recent work on this topic is the generalized adjustment criterion of \cite{perkovic15_uai,perkovic17,perkovic16} which is sufficient but not necessary for the identification of causal effects.  
   \cite{perkovic15_uai, perkovic16,perkovic17} build on prior work of \cite{pearl1993bayesian,shpitser2012validity,vanconstructing} and \cite{maathuis2013generalized}.

  One criterion that is necessary and sufficient for identifying causal effects in DAGs is the g-formula of \cite{robins1986new}. The g-formula is one of the causal identification methods that has seen considerable use in practice  \citep[see e.g.][]{taubman2009intervening, young2011comparative, westreich2012parametric}.  However, the g-formula has not yet been generalized to other types of \mpdag{}s (including CPDAGs).

In this paper, we develop a necessary and sufficient graphical criterion for identifying causal effects  in \mpdag{}s. We refer to our identification criterion  (Theorem \ref{thm:id-formula}) as the causal identification formula. 
The causal identification formula is a generalization of the  g-formula of \cite{robins1986new} to \mpdag{}s. Consequently, we also obtain a generalization of the truncated factorization formula \citep{Pearl2009}, i.e.\ the manipulated density formula \citep{spirtes2000causation} in Corollary \ref{cor:truncfac}. 

From a theoretical perspective, it is of interest to note that the proof of our causal identification formula does not consider intervening on additional variables in the graph (Section \ref{sec:id-proof}). This alleviates concerns of whether such additional interventions are reasonable to assume as possible (see e.g.\  \citealp{vanderweele2014causal,kohler2018eddie}). 

We compare our result to the  generalized adjustment criterion of \cite{perkovic17}  in Section \ref{sec:adjustment}.
Even though the  generalized adjustment criterion is not complete for causal identification, we characterize a special case in which it is ``almost'' complete  in Proposition \ref{prop:adjustment}.

Lastly, \cite{pmlr-v97-jaber19a} recently constructed a graphical algorithm that is necessary and sufficient  for identifying causal effects from observational data that allows for hidden confounders.   The class of graphs that \cite{pmlr-v97-jaber19a} consider is  fully characterized by conditional independences in the observed probability distribution of the data. Their algorithm builds on the work of \cite{tian2002general,shpitser2006identification,huang2006pearl} and \cite{richardson2017nested}.  
 To put our work into wider context, we compare our approach to the approach taken by \cite{pmlr-v97-jaber19a} in  the discussion. 
 Omitted proofs can be found in the Supplement.

\section{PRELIMINARIES} \label{sec:prelims}

We use capital letters (e.g.\ $X$) to denote nodes in a graph as well as random variables that these nodes represent. Similarly, bold capital letters (e.g.\ $\mb{X}$) are used to denote both sets of nodes in a graph as well as the random vectors that these nodes represent.  

\textbf{Nodes, Edges And Subgraphs.} A graph $\g= (\vars,\e) $ consists of a set of nodes (variables) $ \vars=\left\lbrace X_{1},\dots,X_{p}\right\rbrace$ and a set of edges $ \e $. 
The graphs we consider are allowed to contain directed ($\rightarrow$) and  undirected ($-$) edges and at most one edge between any two nodes. An \textit{induced subgraph} $\g_{\mb{V'}} =(\mathbf{V'}, \mathbf{E'})$ of $\g= (\vars,\e) $ consists of  $\mathbf{V'} \subseteq \mathbf{V}$ and  $\mathbf{E'} \subseteq \mathbf{E}$ where $\mathbf{E'}$ are all edges in $\mathbf{E}$ between nodes in $\mathbf{V'}$. An \textit{undirected subgraph} $\g_{undir} = (\mb{V}, \mb{E'})$ of $\g= (\mb{V,E})$ consists of $\mb{V}$ and $\mb{E'} \subseteq \mb{E}$ where $\mb{E'}$ are all undirected edges in $\mb{E}$.

\textbf{Paths.} A \textit{path} $p$ from $X$ to $Y$ in $\g$ is a sequence of distinct nodes $\langle X, \dots,Y \rangle$ in which every pair of successive nodes is adjacent.  A path consisting of undirected edges in an \textit{undirected} path.
A  \textit{causal path} from $X$ to $Y$ is a path from $X$ to $Y$ in which all edges are directed towards $Y$, that is, $X \to \dots\to Y$. Let $p = \langle X=V_0, \dots , V_k=Y \rangle$, $k \geq 1$ be a path in $\g$, $p$ is a \textit{possibly causal path} if no edge $V_{i} \leftarrow V_{j}, 0 \le i < j \le k$ is in $\g$. Otherwise, $p$ is a \textit{non-causal path} in $\g$ (see Definition 3.1 and Lemma 3.2 of \citealp{perkovic17}) (Lemma \ref{lemma:poss-dir-path} in the Supplement). 
For two disjoint subsets $\mathbf{X}$ and $\mathbf{Y}$ of $\mathbf{V}$, a path from $\mathbf{X}$ to $\mathbf{Y}$
is a path from some $X \in \mathbf{X}$ to some $Y \in \mathbf{Y}$.
A path from $\mathbf{X}$ to $\mathbf{Y}$ is \textit{proper} (w.r.t. $\mathbf{X}$) if only its first node is in $\mathbf{X}$.  

\textbf{Partially Directed And Directed Cycles.} 
A causal path from $X$ to $Y$ and the  edge $Y\to X$ form a \textit{directed cycle}.  
A \textit{partially directed cycle} is formed by a possibly causal path from $X$ to $Y$, together with $Y \to X$.

\textbf{Ancestral Relationships.} If $X\to Y$, then $X$ is a \textit{parent} of $Y$. If there is a causal path from $X$ to $Y$, then $X$ is an \textit{ancestor} of $Y$, and $Y$ is a \textit{descendant} of $X$.  If there is a possibly causal path from $X$ to $Y$, then $X$ is a \textit{possible ancestor} of $Y$. We use the convention that every node is a descendant, ancestor, and  possible ancestor of itself.
The sets of parents, ancestors, possible ancestors and descendants of $X$ in~$\g$ are denoted by $\Pa(X,\g)$, $\An(X,\g)$, $\PossAn(X,\g)$ and $\De(X,\g)$ respectively. For a set of nodes $\mathbf{X} \subseteq \mathbf{V}$, we let $\Pa(\mathbf{X},\g) =( \cup_{X \in \mathbf{X}}  \Pa(X,\g)) \setminus \mb{X}$, $\An(\mathbf{X},\g) = \cup_{X \in \mathbf{X}}  \An(X,\g)$, $\PossAn(\mathbf{X},\g) = \cup_{X \in \mathbf{X}}  \PossAn(X,\g)$, and $\De(\mathbf{X},\g) = \cup_{X \in \mathbf{X}}  \De(X,\g)$.

\textbf{Undirected Connected Set.} A node set $\mb{X}$  is an \textit{undirected connected set} in graph $\g$  if for every two distinct nodes $X_i$ and $X_j$ in $\mb{X}$, there is an undirected path from $X_i$ to $X_j$ in $\g$. 

\textbf{Colliders, Shields, And Definite Status Paths.} If a path $p$ contains $X_i \rightarrow X_j \leftarrow X_k$ as a subpath, then $X_j$ is a \textit{collider} on $p$. A path $\langle X_{i},X_{j},X_{k} \rangle$ is an \emph{unshielded triple} if $ X_{i} $ and $ X_{k}$ are not adjacent. A path is \textit{unshielded} if all successive triples on the path are unshielded. A node $X_{j}$ is a \textit{definite non-collider} on a path $p$ if the  edge $X_i \leftarrow X_j$, or the edge $X_j \rightarrow X_k$ is on $p$, or if $X_{i} - X_j - X_k$ is a subpath of $p$ and $X_i$ is not adjacent to $X_k$. A node is of \textit{definite status} on a path if it is a collider, a definite non-collider or an endpoint on the path. A path $p$ is of definite status if every node on $p$ is of definite status.

\textbf{D-connection And Blocking.} A definite status path \textit{p} from $X$ to $Y$ is \textit{d-connecting} given a node set $\mathbf{Z}$ ($X,Y \notin \mathbf{Z}$) if every definite non-collider on $p$ is not in $\mathbf{Z}$, and every collider on $p$ has a descendant in $\mathbf{Z}$. Otherwise, $\mathbf{Z}$ \textit{blocks} $p$. If $\mb{Z}$ blocks all definite status paths between $\mb{X}$ and $\mb{Y}$ in \mpdag{} $\g$, then $\mb{X}$ is \textit{d-separated} from $\mb{Y}$ given $\mb{Z}$ in $\g$ \citep[Lemma C.1 of][]{henckel}.

\textbf{DAGs, PDAGs.}
A \textit{directed graph} contains only directed edges.
A \textit{partially directed graph} may contain both directed and undirected edges.
A directed graph without directed cycles is a \textit{directed acyclic graph $(\DAG)$}.  A \textit{partially directed acyclic graph $(\PDAG)$} is a partially directed graph without directed cycles.

	\begin{figure}[t]
	
	\centering
\begin{subfigure}{.11\textwidth}
		\centering
		\captionsetup[subfigure]{width=80pt}
		\begin{tikzpicture}[>=stealth',shorten >=1pt,auto,node distance=2cm,main node/.style={minimum size=0.6cm,font=\sffamily\Large\bfseries},scale=1.3,transform shape=2]
				\tikzstyle{state}=[inner sep=0.5pt, minimum size=5pt]
		\node[state] (Xia) at (0,0) {};
		\node[state] (Xka) at (0,1.2) {};
		\node[state] (Xja) at (1.2,0) {};
		
		\path[->] (Xka) edge (Xia);
		\draw[-]
		(Xia) edge (Xja);
		\end{tikzpicture}
		\caption{}
\end{subfigure}
\begin{subfigure}{.11\textwidth}
		\centering
		\begin{tikzpicture}[>=stealth',shorten >=1pt,auto,node distance=2cm,main node/.style={minimum size=0.6cm,font=\sffamily\Large\bfseries},scale=1.3,transform shape=2]
				\tikzstyle{state}=[inner sep=0.5pt, minimum size=5pt]
		
		\node[state] (Xic) at (2.8,0) {};
		\node[state] (Xkc) at (2.8,1.2) {};
		\node[state] (Xjc) at (4,0) {};
		
		\path[->] (Xic) edge (Xkc)
		(Xkc) edge (Xjc);
		\draw[-]
		(Xic) edge (Xjc);
		
		\end{tikzpicture}
		\caption{}
\end{subfigure}
\begin{subfigure}{.11\textwidth}
		\centering
		\begin{tikzpicture}[>=stealth',shorten >=1pt,auto,node distance=2cm,main node/.style={minimum size=0.6cm,font=\sffamily\Large\bfseries},scale=1.3,transform shape=2]
				\tikzstyle{state}=[inner sep=0.5pt, minimum size=5pt]

		\node[state] (Xie) at (0,-1) {};
		\node[state] (Xke) at (1.2,-1) {};
		\node[state] (Xle) at (0,-2.2) {};
		\node[state] (Xje) at (1.2,-2.2) {};
		
		\path[->] (Xke) edge (Xje)
		(Xle) edge (Xje);
		
		\draw[-]  (Xie) edge (Xje);
		\path[-]
		(Xke) edge (Xie)
		(Xle) edge (Xie);
		\end{tikzpicture}
		\caption{}
\end{subfigure}
\begin{subfigure}{.11\textwidth}
		\centering
		\begin{tikzpicture}[>=stealth',shorten >=1pt,auto,node distance=2cm,main node/.style={minimum size=0.6cm,font=\sffamily\Large\bfseries},scale=1.3,transform shape=2]
		\tikzstyle{state}=[inner sep=0.5pt, minimum size=5pt]
		
		\node[state] (Xig) at (2.8,-1) {};
		\node[state] (Xjg) at (4,-1) {};
		\node[state] (Xkg) at (2.8,-2.2) {};
		\node[state] (Xlg) at (4,-2.2) {};
		
		\draw[->] (Xlg) edge (Xkg);
	\draw[->] (Xjg) edge (Xlg);
		\draw[-]
		(Xig) edge (Xkg);
		\draw[-]
		(Xig) edge (Xjg);
				\draw[-]
		(Xig) edge (Xlg);
		\end{tikzpicture}
				\caption{}
\end{subfigure}
	\caption{Forbidden induced subgraphs of an MPDAG \citep[see orientation rules in][]{meek1995causal}.}
	\label{fig:orientationRules}
\end{figure}

\noindent\textbf{Markov Equivalence And CPDAGs.}  \citep[c.f.][]{meek1995causal,andersson1997characterization} All DAGs that encode the same d-separation relationships are \textit{Markov equivalent} and form a \textit{Markov equivalence class} of DAGs, which can be \textit{represented} by a \textit{completed partially directed acyclic graph} (CPDAG).

\noindent\textbf{MPDAGs.} A PDAG $\g$ is a \textit{maximally oriented} PDAG (\mpdag{}) if and only if the graphs in Figure \ref{fig:orientationRules} are \textbf{not} induced subgraphs of $\g$. Both a DAG and a CPDAG are types of \mpdag{} \citep{meek1995causal}.

\textbf{$\g$ And $[\g]$.} A DAG $\g[D]$ is \textit{represented} by \mpdag{} $\g$ if $\g[D]$ and $\g$ have the same adjacencies, same unshielded colliders and if for every directed edge $X \rightarrow Y$ in $\g$, $X \rightarrow Y$ is in $\g[D]$ \citep{meek1995causal}.
If $\g$ is an \mpdag{}, then $[\g]$ denotes the set of all $\DAG$s represented by $\g$. 

\textbf{Partial Causal Ordering.} Let $\g[D] = (\mb{V,E})$ be a DAG. A total ordering, $<$, of nodes $\mb{V'} \subseteq \mb{V}$ is \textit{consistent} with $\g[D]$ and called a \textit{causal ordering} of $\mb{V'}$ if for every $X_i, X_j \in \mb{V'}$, such that $X_i < X_j$ and such that $X_i$ and $X_j$ are adjacent in $\g[D]$, $X_i \rightarrow X_j$ is in $\g[D]$. There can be more than one causal ordering of $\mb{V'}$ in a DAG $\g[D] = (\mb{V,E})$. For example, in DAG $X_i \leftarrow X_j \rightarrow X_k$ both orderings $X_j < X_i < X_k$ and $X_j < X_k < X_i$ are consistent.

Let $\g = (\mb{V,E})$ be an \mpdag{}. Since $\g$ may contain undirected edges, there is generally no causal ordering of $\mb{V'}$, for a node set $\mb{V'} \subseteq \mb{V}$ in $\g = (\mb{V,E})$. Instead, we define a \emph{partial causal ordering,} $<$, of $\mb{V'}$ in $\g$ as a total ordering of pairwise disjoint node sets $\mb{A_1}, \dots , \mb{A_k},$ $k \ge 1$, $\cup_{i=1}^{k} \mb{A_i} = \mb{V'}$, that satisfies the following: if $\mb{A_i} < \mb{A_j}$ and there is an edge between $A_i \in \mb{A_i}$ and $A_j \in \mb{A_j}$ in $\g$, then $A_i \rightarrow A_j$ is in $\g$.

\textbf{Do-intervention.} We consider interventions $do(\mathbf{X} =\mathbf{x})$ (for $\mathbf{X}\subseteq \mathbf{V}$) or $do(\mathbf{x})$ for shorthand, which represent outside interventions that set $\mathbf{X}$ to $\mathbf{x}$.

\textbf{Observational And Interventional Densities.} A density $f$ of $\mathbf{V}$ is \textit{consistent} with a $\DAG$ $\g[D] =(\mathbf{V},\mathbf{E})$ if it factorizes as $f(\mb{v})= \prod_{V_i \in \mb{V}}f(v_{i}|\pa(v_{i},\g[D]))$ \citep{Pearl2009}. A density $f$  that is consistent with  $\g[D] =(\mathbf{V},\mathbf{E})$ is also called an \textit{observational density}. 

Let  $\mb{X}$ be a subset of $\mb{V}$ and $\mb{V'}  = \mb{V} \setminus \mb{X}$ in a $\DAG$ $\g[D]$.
A density over $\mb{V'}$ is denoted by $f(\mb{v'}|do(\mb{x}))$, or $f_{\mb{x}}(\mb{v'})$, and called an  \textit{interventional density consistent with}  $\g[D]$ if there is an observational density $f$ consistent with  $\g[D]$ such that $f(\mb{v'}| do(\mb{x}))$  factorizes as
\begin{align}
f(\mathbf{v'}|do(\mathbf{x})) =\prod_{V_{i} \in \mb{V'}}f(v_{i}|\pa(v_{i},\g[D])),
\label{eq11}
\end{align}
for values $\pa(v_i,\g[D])$ of $\Pa(V_i,\g[D])$ that are in agreement with $\mb{x}$.
If $\mb{X} = \emptyset$, we define  $f(\mb{v} | do(\emptyset)) = f(\mb{v})$.
Equation \eqref{eq11} is known as the truncated factorization formula \citep{Pearl2009}, manipulated density formula \citep{spirtes2000causation} or the g-formula \citep{robins1986new}. A density $f$ of $\mb{V}$ is consistent with an \mpdag{} $\g = (\mb{V,E})$ if $f$ is consistent with a DAG in $[\g]$.

 A  density $f(\mb{v'} |do(\mb{x}))$ of $\mb{V'} \subset \mb{V}$, $\mb{X} = \mb{V} \setminus \mb{V'}$ is an \textit{interventional density} consistent with an \mpdag{} $\g = (\mb{V,E})$ if it is an interventional density consistent with a DAG in $[\g]$.
Let $\mb{Y} \subset \mb{V'}$,  and let $f(\mb{v'} |do(\mb{x}))$ be an interventional density consistent with  an \mpdag{}  $\g = (\mb{V,E})$ for some $\mb{X} \subset \mb{V}$, $\mb{V'}  = \mb{V} \setminus \mb{X}$, then $f(\mb{y} |do(\mb{x}))$ denotes the marginal density of $\mb{Y}$ calculated from $f(\mb{v'} |do(\mb{x}))$.

\textbf{Probabilistic Implications Of D-separation.} Let $f$ be any density over $\mb{V}$ consistent with an \mpdag{} $\g = (\mb{V},\mb{E})$ and let $\mb{X,Y},$ and $\mb{Z}$ be pairwise disjoint node sets in $\mb{V}$. If $\mb{X}$ and $\mb{Y}$ are d-separated given $\mb{Z}$ in $\g$, then $\mb{X}$  and $\mb{Y}$ are  conditionally independent  given $\mb{Z}$ in the observational probability density $f$ consistent with $\g[D]$ \citep{lauritzen1990independence,Pearl2009}.

\section{RESULTS} \label{sec:id}

The causal effect of a set of treatments $\mb{X}$ on a set of responses $\mb{Y}$ is a function of the interventional density $f(\mb{y} |do(\mb{x}))$.  
For example, under the assumption of a Bernoulli distributed treatment variable $X$, the causal effect of $X$ on a singleton response $Y$ may be defined as the difference in  expectation of $Y$  under $do(X = 1)$ and $do(X=0)$, that is, $E[Y|do(X =1)] - E[Y|do(X=0)]$ \citep[Chapter 1 in ][]{hernanrobins}. 

We consider a causal effect to be identifiable in an \mpdag{} $\g$ if the interventional density of the response can be uniquely computed from $\g$.
A precise definition is given in Definition \ref{def:causalID}. Definition \ref{def:causalID} is analogous to the Definition 3 of \cite{galles1995testing} and Definition 1 of \cite{pmlr-v97-jaber19a}.

\begin{definition}[\textbf{Identifiability of Causal Effects}] \label{def:causalID} 
Let $\mb{X}$ and $\mb{Y}$ be disjoint node sets in an \mpdag{} $\g = ( \mb{V,E})$. The causal effect of $\mb{X}$  on $\mb{Y}$ is  identifiable in $\g$ if  $f(\mb{y}|do(\mb{x}))$ is uniquely computable from any  observational density consistent with $\g$. \\
Hence, there are no  two DAGs $\g[D]^{1}$, $\g[D]^{2}$ in $[\g]$ such that  
\begin{enumerate}
\item $f_1(\mb{v}) = f_{2}(\mb{v}) = f(\mb{v})$, where $f$ is an observational density consistent with $\g$, and
\item $f_1(\mb{y} | do(\mb{x})) \neq f_{2}(\mb{y}|do(\mb{x}))$, where $f_1(\cdot|do(\mb{x}))$ and $f_2(\cdot|do(\mb{x}))$ are  interventional densities consistent with $\g[D]^1$ and  $\g[D]^2$ respectively.
\end{enumerate}
\end{definition}

\begin{figure}[!tbp]
   \begin{subfigure}{.2\textwidth}
     \centering
     \begin{tikzpicture}[>=stealth',shorten >=1pt,auto,node distance=2cm,main node/.style={minimum size=0.6cm,font=\sffamily\Large\bfseries},scale=0.6,transform shape=2]
   \node[main node]         (X)                        {$X$};
   \node[main node]         (Y) [right of= X]  		{$Y$};
    \draw[-] (X) edge    (Y);
   \end{tikzpicture}
     \caption{}
     \label{proper-cpdag}
   \end{subfigure}
   \unskip
   \vrule
   \begin{subfigure}{.2\textwidth}
     \centering
     \begin{tikzpicture}[>=stealth',shorten >=1pt,auto,node distance=2cm,main node/.style={minimum size=0.6cm,font=\sffamily\Large\bfseries},scale=.6,transform shape=2]
   \node[main node]         (X)                        {$X_1$};
   \node[main node]         (A) [right of= X]  		{$X_2$};	
   \node[main node]       	 (Y)  [right of= A] 	{$Y$};
   \draw[-] (X) edge    (A);
   \draw[->] (A) edge    (Y);
   \end{tikzpicture}
     \caption{}
     \label{proper-mpdag} 
   \end{subfigure}  
    \caption{\small (a) \mpdag{} $\g[C]$, (b) \mpdag{} $\g$.}
   \label{figex2}
 \end{figure}

\subsection{A NECESSARY CONDITION FOR IDENTIFICATION} \label{sec:id-necess}

Proposition \ref{lemmaHedgePDAG} presents a necessary condition for the identifiability of causal effects in \mpdag{}s. This necessary condition is referred to as amenability by \cite{perkovic15_uai,perkovic17}.

\begin{proposition}\label{lemmaHedgePDAG}
 Let $\mb{X}$ and $\mb{Y}$ be disjoint node sets in an \mpdag{} $\g = ( \mb{V,E})$. If there is a proper possibly causal path from $\mb{X}$ to $\mb{Y}$ that starts with an undirected edge in $\g$, then the causal effect of $\mb{X}$  on $\mb{Y}$ is not identifiable in $\g$.
\end{proposition}

Consider MPDAG $\g[C]$ in Figure \ref{proper-cpdag}. Since $X - Y$ is in $\g[C]$, by Proposition \ref{lemmaHedgePDAG}, the causal effect of $X$ on $Y$ is not identifiable in $\g[C]$.  This is intuitively clear since both $X \rightarrow Y$ and $X \leftarrow Y$ are DAGs represented by $\g[C]$. The DAG $X \leftarrow Y$ implies that there is no causal effect of $X$ on $Y$.  
Conversely, the DAG $X \rightarrow Y$  implies that there is a causal effect of $X$ on $Y.$ 

The condition in Proposition \ref{lemmaHedgePDAG} is somewhat less intuitive for non-singleton $\mb{X}$. Consider \mpdag{} $\g$ in Figure \ref{proper-mpdag} and let $\mb{X} = \{X_1, X_2\}$ and $\mb{Y} =  \{Y\}$. The path $X_1 - X_2 \rightarrow Y$ in $\g$ is a possibly causal path from $X_1$ to $Y$ that starts with an undirected edge. However,   $X_1 - X_2 \rightarrow Y$ is not a proper possibly causal path from $\mb{X}$ to $Y$, since it contains $X_2$ and $X_1$. Hence, the causal effect of $\mb{X}$ on $Y$ may still be identifiable in $\g$.

\subsection{PARTIAL CAUSAL ORDERING IN MPDAGS} \label{sec:id-pto}

For the proof of our main result, it is necessary to determine a partial causal ordering for a set of nodes in an \mpdag{}.
In order to compute a partial causal ordering of nodes in an \mpdag{}, we first define a bucket.

 \begin{definition}[\textbf{Bucket}] \label{def:max-bucket} 
Let $\mb{D}$ be a node set in an \mpdag{} $\g = ( \mb{V,E})$. If $\mb{B}$ is a maximal undirected connected subset of $\mb{D}$ in $\g$, we call $\mb{B}$ a \textit{bucket} in $\mb{D}$.
\end{definition}

 Definition \ref{def:max-bucket} is similar to the definition of a bucket of  \cite{jaber2018causal}.
One difference 
is that Definition \ref{def:max-bucket} allows for directed edges between the nodes within the same bucket, whereas the definition of \cite{jaber2018causal} does not. For instance, $\{X,V_1,Y_1\}$ is a bucket in $\mb{V}$ in MPDAG $\g = (\mb{V,E})$ in Figure \ref{fig:example_meekgraph11}(a). Note that since we require a bucket to be a maximal undirected connected set, $\{X,V_1\}$ is not a bucket in $\mb{V}$.
 
 Definition \ref{def:max-bucket}  can be used to induce a unique partition of any  node set $\mb{D}$ in an MPDAG $\g = (\mb{V,E})$, $\mb{D} \subseteq	 \mb{V}$. We refer to this partition as \textit{the bucket decomposition} in the corollary of Definition \ref{def:max-bucket} below.
 
 \begin{corollary}[\textbf{Bucket Decomposition}]
Let $\mb{D}$ be a node set in an \mpdag{} $\g = ( \mb{V,E})$. Then there is a unique partition of $\mb{D}$ into $\mb{B_1}, \dots \mb{B_k} $, $k \ge 1$  in $\g$ induced by Definition \ref{def:max-bucket}. That is 
\vspace{-.3cm}
\begin{itemize}[leftmargin=*]
\item $\mb{D} = \cup_{i=1}^k \mb{B_i}$, and
\item  $\mb{B_i} \cap \mb{B_j} = \emptyset$,  $i,j \in \{ 1, \dots , k\}$,  $i \neq j$, and 
\item  $\mb{B_i}$ is a bucket in $\mb{D}$ for each $i \in \{ 1, \dots , k\}$.
\end{itemize}
\label{cor:bucket}
\end{corollary}

\begin{figure}
\tikzstyle{every edge}=[draw,>=stealth',->]
\newcommand\dagvariant[1]{\begin{tikzpicture}[xscale=.5,yscale=0.5]
\node (a) at (0,0) {};
\node (d) at (0,2) {};
\node (b) at (2,0) {};
\node (c) at (2,2) {};
\begin{scope}[gray]
\draw (a) edge [-] (b);
\draw (b) edge [-] (c);
\draw (d) edge [-] (a);
\draw (d) edge [-] (b);
\draw (a) edge [-] (c);
\end{scope}
\draw #1;
\end{tikzpicture}}

\begin{subfigure}{.4\columnwidth}
  \centering
\begin{tikzpicture}[->,>=latex,shorten >=1pt,auto,node distance=0.8cm,scale=1,transform shape]
  \tikzstyle{state}=[inner sep=1pt, minimum size=12pt]

  \node[state] (Xia) at (0,0) {$Y_1$};
  \node[state] (Xka) at (0,2) {$V_1$};
  \node[state] (Xja) at (2,0) {$X$};
  \node[state] (Xsa) at (2,2) {$Y_2$};

  \draw (Xka) edge [-] (Xja);
\draw 		(Xia) edge  [-] (Xka);
 \draw   	(Xja) edge [->,line width=1.1pt] (Xsa);
 \draw    	(Xia) edge [->,line width=1.1pt] (Xja);
\draw    	(Xia) edge [->,line width=1.1pt] (Xsa);
\end{tikzpicture}
\caption{}
  \label{mpdag}
\end{subfigure}
\unskip
\vrule
\hspace{0.5cm}
\begin{subfigure}{.5\columnwidth}
\dagvariant{
(a)  edge [->,line width=1.3pt,color=blue]  (b)
(b)  edge [->,line width=1.3pt,color=blue]  (c)
(a) edge [->,line width=1.3pt,color=blue] (c)
(d) edge [->] (a)
(d) edge [->] (b)
}
\dagvariant{
(a)  edge [->,line width=1.3pt,color=blue]  (b)
(b)  edge [->,line width=1.3pt,color=blue]  (c)
(a) edge [->,line width=1.3pt,color=blue] (c)
(a) edge [->] (d)
(d) edge [->] (b)
}
\dagvariant{
(a)  edge [->,line width=1.3pt,color=blue]  (b)
(b)  edge [->,line width=1.3pt,color=blue]  (c)
(a) edge [->,line width=1.3pt,color=blue] (c)
(a) edge [->] (d)
(b) edge [->] (d)
}
\caption{}
\label{alldags-mpdag}
\end{subfigure}
\caption{(a) $\MPDAG$ $\g$, (b) $\DAG$s represented by $\g$.} 
\label{fig:example_meekgraph11}
\end{figure}

Consider \mpdag{} $\g = (\mb{V,E})$ in Figure \ref{mpdag}.  In order to find the bucket decomposition of $\mb{V}$ in $\g$, let us consider  the  undirected subgraph $\g_{undir}$ of $\g$. The only path in $\g_{undir}$ is   $Y_1-V_1 -X$. Hence, the bucket decomposition of $\mb{V}$ is  $ \{ \{X,V_1, Y_1\}, \{Y_2\}\}$.  

Consider DAGs in Figure \ref{alldags-mpdag}, which are all DAGs represented by $\g$  in Figure \ref{mpdag}.  Some total orderings of $\mb{V}$ that are consistent with DAGs in Figure   \ref{alldags-mpdag} are: $V_1 <Y_1  < X< Y_2$, $Y_1 < V_1 < X < Y_2 $, and $Y_1 < X < V_1 < Y_2$, from left to right respectively. These three orderings are consistent with the following partial causal ordering $\{X,V_1,Y_1\} < Y_2$, which is a total ordering of the buckets in the bucket decomposition of $\mb{V}$. This motivates Algorithm \ref{alg:pto}.

Algorithm \ref{alg:pto} outputs an ordered bucket decomposition of a set of nodes $\mb{D}$ in an \mpdag{} $\g$. The proof that Algorithm \ref{alg:pto} will always complete is given in Lemma \ref{lemma:algo-completes} in the Supplement. 
Next, we prove that the ordered list of buckets  output by Algorithm \ref{alg:pto} is a partial causal ordering of $\mb{D}$ in $\g$ (Lemma \ref{lemma:PCO}).  Algorithm \ref{alg:pto} and Lemma \ref{lemma:PCO} are similar to the PTO algorithm and  Lemma 1 of \cite{jaber2018graphical}.

\begin{lemma} \label{lemma:PCO}
Let $\mb{D}$ be a node set in an \mpdag{} $\g = (\mb{V,E})$ and let $ (\mb{B_1}, \dots, \mb{B_k})$, $k \ge 1$, be the output of \texttt{PCO}$(\mb{D}, \g)$. Then for each $i,j \in \{1, \dots k\}$, $\mb{B_i}$ and $\mb{B_j}$ are buckets in $\mb{D}$ and if $i < j$, then $\mb{B_i} < \mb{B_j}$ in $\g$. \end{lemma}

\begin{algorithm}[t]
\vspace{.2cm}
 \Input{Node set $\mb{D}$ in  \mpdag{} $\mathcal{G}{=}(\mathbf{V,E})$.}
 \Output{An ordered list $\mb{B} {=} (\mb{B_1}, \dots , \mb{B_k}), k \ge 1,$ of the bucket decomposition of $\mb{D}$ in $\g$.}
 \vspace{.2cm}
   \SetAlgoLined
   Let $\mb{ConComp}$ be the bucket decomposition of $\mb{V}$ in $\g$;\\
   Let $\mb{B}$ be an empty list;\\
\While{$\mb{ConComp}\neq \emptyset$}{
Let  $\mb{C} \in \mb{ConComp}$;\\
Let $\overline{\mb{C}}$ be the set of nodes in $\mb{ConComp}$ that are not in $\mb{C}$;\\
\If{all edges between  $\mb{C}$ and  $\overline{\mb{C}}$   are into $\mb{C}$ in $\g$}
{Remove $\mb{C}$ from $\mb{ConComp}$;\\
Let $\mb{B_{*}} = \mb{C} \cap \mb{D}$;\\
\If{$\mb{B_{*}}\neq \emptyset$} { Add $\mb{B_{*}}$ to the beginning of $\mb{B}$;}
}
}
\Return  $\mb{B}$;\\
\caption{Partial causal ordering (PCO)}
\label{alg:pto}
\end{algorithm}

Consider  \mpdag{} $\g = (\mb{V,E})$ in Figure \ref{mpdag} and let $\mb{D} = \{X,Y_1,Y_2\}$.  We now explain how the output of \texttt{PCO}$(\mb{D},\g)$ is obtained. 

In line 2, the bucket decomposition of $\mb{V}$ is obtained, $\mb{ConComp} = \{ \{X,Y_1,V_1\}, \{Y_2\}\}$ (as noted above). In line 2, $\mb{B}$ is initialized as an empty list.

Let $\mb{C} =  \{X,Y_1,V_1\}$ (line 4). Then $\overline{\mb{C}} = \{Y_2\}$ (line 5). Since $X \rightarrow Y_2$ and $Y_1 \rightarrow Y_2$ are in $\g$, $\mb{C}$ does not satisfy the condition in line 6 and hence, $ \{X,Y_1,V_1\}$ cannot be removed from $\mb{ConComp}$ at this time. 

Next, $\mb{C} = \{Y_2\}$ (line 4) and  $\overline{\mb{C}} = \{X,Y_1,V_1\}$ (line 5).
 Since all edges between $\{Y_2\}$ and $\{X,Y_1,V_1\}$  in $\g$ are into $\{Y_2\}$, Algorithm \ref{alg:pto} removes $\{Y_2\}$ from $\mb{ConComp}$ in line 7. Since $\mb{B_{*}} = \mb{C} \cap \mb{D} = \{Y_2\}$ (line 8), Algorithm \ref{alg:pto} adds $\{Y_2\}$ to the beginning of list $\mb{B}$ (line 10).

Now, $\mb{C} =  \{X,Y_1,V_1\}$  (line 4) and $\overline{\mb{C}} = \emptyset$ (line 5).  Hence, $\mb{C}$ satisfies condition in line 6 and  $\mb{C}$ is removed from $\mb{ConComp}$ (line 7).  Then $\mb{B_{*}} = \mb{C} \cap \mb{D} = \{X,Y_1\}$ (line 8), and $\mb{B} = (\{X,Y_1\}, \{Y_2\})$ (line 10). 
Since $\mb{ConComp}$ is empty, Algorithm \ref{alg:pto} outputs $\mb{B}$.

\subsection{CAUSAL IDENTIFICATION FORMULA} \label{sec:factorid}

We present our main result  which we refer to as the causal identification formula in Theorem \ref{thm:id-formula}. Theorem \ref{thm:id-formula} establishes that the condition from Proposition \ref{lemmaHedgePDAG} is not only necessary, but also sufficient for the identification of causal effects in \mpdag{}s.

\begin{theorem}[\textbf{Causal identification formula}]
\label{thm:id-formula}
Let $\mb{X}$ and $\mb{Y}$ be disjoint node sets in an  \mpdag{} $\g = (\mb{V},\mb{E})$. 
If   there is no proper possibly causal path from $\mb{X}$ to $\mb{Y}$ in $\g$ that starts with an undirected edge, then for any  observational density $f$  consistent with $\g$ we have
\begin{align} \label{gform}
f(\mb{y}|do(\mb{x})) = \int \prod_{i=1}^{k} f(\mb{b_i} | \pa(\mb{b_i},\g))d\mb{b},
\end{align}
for values $\pa(\mb{b_i},\g)$ of $\Pa(\mb{b_i},\g)$ that are in agreement with $\mb{x}$,
where $(\mb{B_1}, \dots ,\mb{B_k})=$ \texttt{PCO}$(\An(\mb{Y}, \g_{\mb{V} \setminus \mb{X}}), \g)$ and $\mb{B} =\An(\mb{Y}, \g_{\mb{V} \setminus \mb{X}}) \setminus \mb{Y}$.
\end{theorem}

For a DAG $\g[D] = (\mb{V,E})$, it is well known that in order to identify a causal effect of $\mb{X}$ on $\mb{Y}$ in $\g[D]$ it is enough to consider the set of ancestors of $\mb{Y}$, that is $\An(\mb{Y}, \g[D])$ (see Theorem 4 of \citealp{tian2002general}).  
The causal identification formula refines this notion  by using a subset of ancestors of $\mb{Y}$ to identify the causal effect of $\mb{X}$ on $\mb{Y}$ in an \mpdag{} $\g$. The variables that appear on the right hand side of equation \eqref{gform} are  in $\An(\mb{Y}, \g_{\mb{V} \setminus \mb{X}})$, or in $\mb{X}$, for those $\mb{X}$ that have a proper causal path to $\mb{Y}$ in $\g$.

  The causal identification formula is a generalization of  the g-formula of \cite{robins1986new}, the truncated factorization formula of \cite{Pearl2009}, or the manipulated density formula of \cite{spirtes2000causation} to the case of \mpdag{}s. To further exhibit this connection, we include the following corollary.  

\begin{corollary}[\textbf{Factorization and truncated factorization formula in \mpdag{}s}] \label{cor:truncfac}
Let $\mb{X}$ be a node set in an \mpdag{} $\g = (\mb{V,E})$ and let $\mb{V'} = \mb{V} \setminus \mb{X}$. Furthermore, let $(\mb{V_1}, \dots ,\mb{V_k})$ be the output of  \texttt{PCO}$(\mb{V} ,\g)$.
Then for any observational density $f$  consistent with $\g$ we have
\begin{enumerate}
\item $ f(\mb{v}) =  \prod_{\mb{V_i} \subseteq \mb{V}} f(\mb{v_i} | \pa(\mb{v_i},\g)),$
\item If there is no pair of nodes $V \in \mb{V'}$ and $X \in \mb{X}$ such that $X - V$ is in $\g$, then 
$$ f(\mb{v'}|do(\mb{x})) =  \prod_{\mb{V_i} \subseteq \mb{V'}} f(\mb{v_i} | \pa(\mb{v_i},\g)),$$
for values $\pa(\mb{v_i},\g)$ of $\Pa(\mb{v_i},\g)$ that are in agreement with $\mb{x}$.
\end{enumerate}
\end{corollary}

Whenever $f(\mb{v'}| do(\mb{x}))$ is identifiable in \mpdag{} $\g = (\mb{V,E})$, we can also identify $f(\mb{y}| do(\mb{x}))$ as
$$f(\mb{y}| do(\mb{x})) = \int f(\mb{v'}|do(\mb{x})) d \overline{\mb{v'}},  $$
where $\mb{X}$ and $\mb{Y}$ are disjoint subsets of $\mb{V}$, $\mb{V'} = \mb{V} \setminus \mb{X}$, and $
\overline{\mb{V'}} = \mb{V} \setminus \{ \mb{X}\cup \mb{Y}\}$.
Since the necessary condition for identifying $f(\mb{v'}| do(\mb{x}))$ (Corollary \ref{cor:truncfac}) is generally stronger than the necessary condition for identifying $f(\mb{y}|do(\mb{x}))$ there are cases when  $f(\mb{y}|do(\mb{x}))$ is identifiable and  $f(\mb{v'}| do(\mb{x}))$ is not identifiable. One such case  is explored in Example \ref{ex:motivate}.

\subsection{EXAMPLES} \label{examples}

\begin{figure}
   \centering
      \begin{subfigure}{.25\textwidth}
  \vspace{1cm}
     \centering
     \begin{tikzpicture}[>=stealth',shorten >=1pt,auto,node distance=2cm,main node/.style={minimum size=0.8cm,font=\sffamily\Large\bfseries},scale=0.75,transform shape]
   \node[main node]         (X)                       {$X$};
   \node[main node]         (V1) [above right of= X]  		{$V_{1}$};
   \node[main node]         (V2) [right of = V1]  		{$V_2$};
      \node[main node]         (V3) [above of = X]  		{$V_3$};
   \node[main node]       	 (Y)  [below right of= V2]        {$Y$};
      \draw[->] (V3) edge    (X);
         \draw[-] (V3) edge    (V1);
   \draw[->] (X) edge    (Y);
   \draw[-] (V1) edge    (V2);
   \draw[->] (V2) edge    (X);
      \draw[->] (V2) edge    (Y);
   \draw[->] (V1) edge    (X);
   \draw[->] (V1) edge    (Y);
   \end{tikzpicture}
    \caption{}
    \label{ex1}
   \end{subfigure}
      \vrule
      \begin{subfigure}{.22\textwidth}
     \centering
     \begin{tikzpicture}[>=stealth',shorten >=1pt,auto,node distance=2cm,main node/.style={minimum size=0.8cm,font=\sffamily\Large\bfseries},scale=.7,transform shape]
   \node[main node]         (X1)                       {$X_1$};
   \node[main node]         (V1) [right of= X1]  		{$V_{4}$};
   \node[main node]         (Y) [right of = V1]  		{$Y$};
   \node[main node]         (X2) [below of= V1]  		{$X_{2}$};
      \node[main node]         (V2) [left of= X2]  		{$V_{3}$};
         \node[main node]         (V3) [above of= X1]  		{$V_{1}$};
  \node[main node]         (V4) at (1,2)		{$V_{2}$};
   \draw[->] (X1) edge    (V1);
   \draw[->] (V1) edge    (Y);
   \draw[->] (V1) edge    (X2);
   \draw[->] (X2) edge    (Y);
   \draw[->] (V4) edge    (X1);
   \draw[->] (V3) edge    (X1);
     \draw[->] (V2) edge    (X2);     
   \draw[->] (X1) ..  controls (1.5,1.7) and (2.5,1.7) ..     (Y);
   \end{tikzpicture}
    \caption{}
    \label{ex2}
   \end{subfigure}
      \caption{\small (a) \Mpdag{} $\g$, (b) $\DAG$ $\g[D]$.}
   \label{fig:nogbc1}
\end{figure}

\begin{example}\label{ex:motivate} In this example,  the causal effect of $\mb{X}$ on $\mb{Y}$ is identifiable in an \mpdag{} $\g = (\mb{V}, \mb{E})$, but the effect of $\mb{X}$ on $\mb{V'} = \mb{V} \setminus \mb{X}$ is not identifiable in $\g$.

Consider \mpdag{} $\g$ in Figure \ref{mpdag}  and let $f$ be an observational density consistent with $\g$. Let  $\mathbf{X} = \{X\}$ and $\mathbf{Y} = \{Y_1,Y_2\}$.  Note that path $X - V_1 - Y_1$ while proper is not possibly causal from $X$ to $Y_1$ in $\g$ due to edge $Y_1 \rightarrow X$.
The only possibly causal path from $X$ to $\mb{Y}$ in $\g$ is $X \rightarrow Y_2$. Hence, by Theorem \ref{thm:id-formula}, the causal effect of $X$ on $\mb{Y}$ is identifiable in $\g$.  

To use the causal identification formula we first determine that $\An(\{Y_1,Y_2\}, \g_{\mb{V} \setminus \{X\}}) = \{Y_1, Y_2\}$, the bucket decomposition of $\{Y_1,Y_2\}$ is $\{\{Y_1\},\{Y_2\}\} $, and  \texttt{PCO}$(\{Y_1, Y_2\},\g) = (\{Y_1\},\{Y_2\})$. 
Next, $\Pa(Y_1 ,\g) = \emptyset$, and $ \Pa(Y_2,\g) = \{X, Y_1\}$.
Hence, by Theorem \ref{thm:id-formula}, 
$ f(y_1,y_2|do(x)) =  f(y_2|x, y_1) f(y_1)$.

Now, let $\mb{V'} = \mb{V} \setminus \{X\}$. Since  $X - V_1$ is in $\g$, by Corollary \ref{cor:truncfac}, $f(\mb{v'} |do(\mb{x}))$ is not identifiable in $\g$. 
\end{example}

\begin{example}  In this example, both  the causal effect of $X$ on $Y$ and the causal effect of $X$ on $\mb{V'} = \mb{V} \setminus \{X\}$ are identifiable in an \mpdag{} $\g = (\mb{V}, \mb{E})$.

Consider \mpdag{} $\g$ in Figure \ref{ex1}  and let $f$ be an observational density consistent with $\g$. The only possibly causal path from $X$ to $Y$ in $\g$ is $X \rightarrow Y$. Hence, the causal effect of $X$ on $Y$ is identifiable in $\g$. 

In fact, there are no undirected edges connected to $X$, so the causal effect of $X$ on $\mb{V'}$, $\mb{V'} = \{V_1, V_2,V_3, Y\}$ is also identifiable in $\g$. Thus, we can obtain the truncated factorization formula with respect to $X$ in $\g$.

We will first determine the causal identification formula for $f(y|do(x))$ in $\g$. We first identify that $\An(Y, \g_{\mb{V} \setminus \{X\}}) = \{V_1, V_2, Y\}$. 
The bucket decomposition of  $\{V_1, V_2, Y\}$ is $\{\{V_1,V_2\},\{Y\}\}$ and \texttt{PCO}$(\{V_1, V_2, Y\},\g)$ is $(\{V_1,V_2\},\{Y\})$. Furthermore, $\Pa(\{V_1,V_2\} ,\g) = \emptyset,$  $\Pa(Y,\g) = \{X, V_1,V_2\}$.
Hence, by Theorem \ref{thm:id-formula}, the causal identification formula for   $f(y|do(x))$ in $\g$ is 
$ f(y|do(x)) = \int f(y|x,v_1,v_2) f(v_1,v_2)dv_1dv_2.$

To use Corollary \ref{cor:truncfac}, first note that the output of \texttt{PCO}$(\mb{V},\g)$ is  $(\{V_1,V_2,V_3\},\{X\}, \{Y\})$ and that the ordered bucket decoposition of $\mb{V'}$ is $(\{V_1,V_2,V_3\}, \{Y\})$.  Further, $\Pa(\{V_1,V_2,V_3\} ,\g) = \emptyset$. Then, $f(\mb{v'}|do(x)) =  f(y|x,v_1,v_2) f(v_1,v_2,v_3).$
   \label{ex:regular}
\end{example}

\begin{example} This example shows how the causal identification formula can be used to estimate the causal effect of  $\mb{X}$ on $\mb{Y}$ in an \mpdag{} $\g$ under the assumption that the observational density $f$ consistent with $\g$ is multivariate Gaussian. 

Consider DAG $\g[D]$ in Figure \ref{ex2} and let $f$ be an observational density consistent with $\g[D]$. Further, let  $\mathbf{X} = \{X_1,X_2\}$ and $\mathbf{Y} = \{Y\}$.
Then $\An(Y , \g[D]_{\mb{V} \setminus \mb{X}}) = \{Y,V_4\}$, the bucket decomposition of $\{Y,V_4\}$ is $\{\{V_4\},\{Y\}\}$,  and  \texttt{PCO}$( \{Y,V_4\},\g[D]) = (\{V_4\},\{Y\})$ in $\g[D]$.

Since $\Pa(V_4, \g[D]) =\{X_1\}$, and $\Pa(Y, \g[D]) = \{ X_1, X_2, V_4\}$, by Theorem \ref{thm:id-formula},
$$f(y|do(x_1, x_2)) = \int f(y|x_1,x_2,v_4)f(v_4|x_1) dv_4.$$
Suppose that the density $f$ consistent with $\g[D]$ is multivariate Gaussian. The causal effect of $\mb{X}$ on $Y$ can then be defined as the vector 
$$\bigg(\frac{\partial E[Y |do(x_1,x_2)]}{\partial x_1},\frac{\partial E[Y |do(x_1,x_2)]}{\partial x_2}\bigg)^T,$$ \citep{nandy2014estimating}.
Hence, consider $E[Y|do(x_1,x_2)]$,
\begin{align*}
&E[Y|do(x_1,x_2)] = \int  y f(y|do(x_1,x_2))dy\\
& =  \int \int  y f(y|x_1,x_2,v_4)f(v_4|x_1) dv_4dy \\
&= \int E[Y|x_1,x_2,v_4]f(v_4|x_1) dv_4 \\
&= \alpha x_1 + \beta x_2 + \gamma \int v_4 f(v_4|x_1)dv_4 \\
& =  \beta  x_2 +   x_1(\alpha  + \gamma \delta),
\end{align*}
where $E[Y|x_1,x_2,v_4] = \alpha x_1 + \beta x_2 + \gamma v_4$ and $E[V_4|x_1]  = \delta x_1$ \citep[Theorem 3.2.4 of][see Theorem \ref{theorem:mardia-condexp} in the Supplement]{mardia1980multivariate}. 

The causal effect of $\mb{X}$ on $\mb{Y}$ is equal to $(\alpha + \gamma \delta, \beta)$. Consistent estimators for $\alpha$, $\beta$, and $\gamma$ are the least-squares estimators of the respective coefficients of $X_1$, $X_2$, and $V_4$ in the regression of  $Y$ on $X_1$, $X_2$, and $V_4$. Analogously, the consistent estimator for $\delta$ is the least-squares coefficient of $X_1$ in the regression of $V_4$ on $X_1$.
   \label{ex:gauss}
\end{example}
\subsection{PROOF OF THEOREM \ref{thm:id-formula}} \label{sec:id-proof}

The proof of Theorem \ref{thm:id-formula} relies on Lemma \ref{lemma:big2} in the Supplement. Lemma \ref{lemma:big2}  is proven through use of do-calculus  \citep{Pearl2009} and basic probability calculus.

The proofs of  Theorem \ref{thm:id-formula} and Lemma \ref{lemma:big2} do not require intervening on additional variables in $\g$. 
This fact alleviates any concerns of whether such additional interventions are reasonable to assume as possible (see e.g. \citealp{vanderweele2014causal,kohler2018eddie}).
\begin{proofof}[Theorem~\ref{thm:id-formula}]
For $i \in \{2, \dots, k\}$, let $\mb{P_i} = (\cup_{j=1}^{i-1}\mb{B_i}) \cap \Pa(\mb{B_i}, \g)$. For $i \in \{1, \dots, k\}$, let $\mb{X_{p_i}} = \mb{X} \cap \Pa(\mb{B_i}, \g)$.

Then
\begin{align}
&f(\mb{y}|do(\mb{x}))  = \int f(\mb{b}, \mb{y}|do(\mb{x})) d\mb{b} \nonumber\\
&= \int f(\mb{b_1}|do (\mb{x})) \prod_{i=2}^{k} f(\mb{b_{i}}| \mb{b_{i-1}}, \dots, \mb{b_1}, do(\mb{x}))d\mb{b} \nonumber
\\
& = \int f(\mb{b_1}|do (\mb{x})) \prod_{i=2}^{k} f(\mb{b_{i}}| \mb{p_i},do(\mb{x}))d\mb{b} \label{eq:dorule22}\\
& = \int f(\mb{b_1}|do (\mb{x_{p_1}})) \prod_{i=2}^{k} f(\mb{b_{i}}| \mb{p_i}, do(\mb{x_{p_i}}))d\mb{b} \label{eq:dorule32}\\
&=  \int \prod_{i=1}^{k} f(\mb{b_{i}}| pa(\mb{b_i},\g))d\mb{b}, \label{eq:doparents2}
\end{align}
The first two equalities follow from the law of total probability and the chain rule. Equations \eqref{eq:dorule22}, \eqref{eq:dorule32}, and  \eqref{eq:doparents2} follow by  applying  results   \ref{lemma:dorule2ID2},  \ref{lemma:dorule3ID2}, and \ref{lemma:dorulepaID2}  in  Lemma \ref{lemma:big2} in the Supplement. 
\end{proofof}

\section{COMPARISON TO ADJUSTMENT} \label{sec:adjustment}

The current state-of-the-art method for identifying causal effects in \mpdag{}s is the generalized adjustment criterion of \cite{perkovic17} stated in Theorem \ref{defadjustmentmpdag}. 

\begin{theorem}[\textbf{Adjustment set, Generalized adjustment criterion}; \citealp{perkovic17}]
   Let $\mathbf{X,Y}$ and $\mathbf{Z}$ be pairwise disjoint node sets in an  \mpdag{} $\g = (\mb{V},\mb{E})$. Let $f$ be any  observational density consistent with $\g$. 
   
   Then $\mb{Z}$ is an adjustment set relative to $(\mb{X,Y})$ in $\g$ and we have
   \begin{equation}
   f(\mathbf{y}|do(\mathbf{x}))=
   \begin{cases}
   \int f(\mathbf{y}|\mathbf{x,z})f(\mathbf{z})d\mathbf{z}  & \text{, if }\mathbf{Z} \neq \emptyset,\\
    f(\mathbf{y}|\mathbf{x}) & \text{, otherwise.}
   \end{cases}
   \nonumber
   \end{equation}
if  and only if the following conditions are satisfied:
   \begin{enumerate}
     \item\label{cond0} There is no proper possibly causal path from $\mb{X}$ to $\mb{Y}$ that starts with an undirected edge in $\g$.
   \item\label{cond1} $\mathbf{Z} \cap \fb{\g} = \emptyset$,  where
   \begin{align}
   &\fb{\g} = \{ W' \in \mb{V}: W' \in \PossDe(W,\g), \notag\\
   &\text{for some } W \notin \mathbf{X} \, \text{which lies on a proper possibly } \notag\\
   & \text{causal path  from } \mathbf{X} \,\text{to}\, \mathbf{Y}\text{in } \g\}.\notag
\end{align}
   \item\label{cond2} All proper non-causal definite status paths from $\mathbf{X}$ to $\mathbf{Y}$ are blocked by $\mathbf{Z}$~in~$\g$.
   \end{enumerate}
   \label{defadjustmentmpdag}
\end{theorem}

The generalized adjustment criterion is sufficient for identifying causal effects in an  \mpdag{}, but it is not necessary. However, when $\mb{X}$ and $\mb{Y}$ are singleton sets, the generalized adjustment criterion identifies all non-zero causal effects of $\mb{X}$ on $\mb{Y}$ in an \mpdag{} $\g$. This is shown in the following proposition.

\begin{proposition}\label{prop:adjustment}
Let $X$ and $Y$ be distinct nodes in an \mpdag{} $\g = (\mb{V},\mb{E})$.  If $Y \notin \Pa(X,\g)$, then the causal effect of $X$ on $Y$ is identifiable in $\g$ if and only if there is an adjustment set relative to $(X,Y)$ in $\g$.

Furthermore, if $Y \notin \Pa(X,\g)$,  then  $\Pa(X, \g)$  is an adjustment set relative to $(X,Y)$ in $\g$ whenever one such set exists.
\end{proposition}

If $Y \in \Pa(X,\g)$, then due to the acyclicity of $\g$,  there is no causal path from $X$ to $Y$ in $\g$ and therefore no causal effect of $X$ on $Y$  (see Lemma \ref{lemma:preproc} in the Supplement). Hence, by Proposition \ref{prop:adjustment}, the generalized adjustment criterion is ``almost'' complete for the identification of causal effects of variable $X$ on a response $Y$ in \mpdag{}s.

If $\mb{X}$, or $\mb{Y}$ are non-singleton sets in $\g$, however, the generalized adjustment criterion will fail to identify some non-zero causal effects of $\mb{X}$ on $\mb{Y}$. We discuss this further in the two examples below.

\begin{example}\label{ex:gacmotivate}
Consider  \mpdag{} $\g$ in Figure \ref{mpdag} and let $\mb{X} = \{X\}$, and $\mb{Y} = \{Y_1, Y_2\}$ as in Example \ref{ex:motivate}.  

Path $X \leftarrow Y_1$ is a non-causal path from $X$ to $\mb{Y}$ that cannot be blocked by any set of nodes disjoint with  $\{X, Y_1\}$.
 Hence,  there is no adjustment set relative to $(X,\mb{Y})$ in $\g$. But there is a causal path from $X$ to $\mb{Y}$ in $\g$ and  as we have seen in Example \ref{ex:motivate}, the causal effect of $X$ on $\mb{Y}$ is identifiable in $\g$.
\end{example}

\begin{example}
Consider DAG $\g[D]$ in Figure \ref{ex2} and let $\mb{X} = \{X_1, X_2\}$, and $\mb{Y} = \{Y\}$. Then $\fb{\g[D]} = \{V_4,Y\}$.
For a set $\mb{Z}$ to satisfy the generalized adjustment criterion relative to $(\mb{X},Y)$ in $\g$, $\mb{Z}$ cannot contain nodes in $\{V_4, Y\}$, or $\{X_1,X_2\}$ and $\mb{Z}$ must block all proper non-causal paths from $\mb{X}$ to $Y$ in $\g[D]$. 

However, $X_2 \leftarrow V_4 \rightarrow Y$ is a proper non-causal path from $\mb{X}$ to $Y$ in $\g[D]$ that cannot be blocked by any set $\mb{Z}$ that satisfies $\mb{Z} \cap \{X_1,X_2,V_4,Y\} = \emptyset$. 
Hence, there is no adjustment set relative to $(\mb{X},Y)$ in $\g[D]$.  But as we have seen in Example \ref{ex:gauss}, the causal effect of $\mb{X}$ on $Y$ is identifiable in $\g[D]$ and furthermore, both $X_1$ and $X_2$ are causes of $Y$ in $\g[D]$.
\end{example}

\section{DISCUSSION} \label{sec:discussion}

We introduced a causal identification formula that allows complete identification of causal effects in \mpdag{}s.  Furthermore, we gave a comparison of our graphical criterion to the current state of the art method for causal identification in \mpdag{}s.

Since the causal identification formula comes in the familiar form of the g-formula of \cite{robins1986new} for DAGs, our results can be used to generalize applications of the g-formula to \mpdag{}s. For example,   \cite{murphy2003optimal},  \cite{collins2004conceptual},  and \cite{collins2007multiphase} give criteria for estimating the optimal dynamic treatment regime from longitudinal data that are based on the g-formula.  This idea can further be combined with recent work of  \cite{rahmadi2017causality} and \cite{rahmadi2018causality} that establishes an approach for estimating the \mpdag{} using data from longitudinal studies.

Throughout the paper, we assume no latent variables.
When latent variables are present,  one can at most  learn a partial ancestral graph (PAG)  over the set of observed variables from the observed data \citep{richardson2002ancestral,spirtes2000causation,zhang2008causal,zhang2008completeness}. PAGs represent an equivalence class of DAGs over the set  of observed and unobserved  variables. 

 \cite{pmlr-v97-jaber19a} recently developed a recursive graphical algorithm that is both necessary and sufficient for identifying causal effects in PAGs. 
 Our causal identification formula does not follow as a simplification of the result of \cite{pmlr-v97-jaber19a}. To see this, notice that the strategy of \cite{pmlr-v97-jaber19a} for identifying causal effects in PAG $\g[P]$ relies on the fact that the causal effect of $\mb{X}$ on $\mb{Y}$  is identifiable in $\g[P]$ if and only if the causal effect of $\mb{V} \setminus \PossAn(\mb{Y}, \g[P]_{\mb{V} \setminus \mb{X}})$ on  $ \PossAn(\mb{Y}, \g[P]_{\mb{V} \setminus \mb{X}})$ is identifiable  in $\g[P]$ (see equation (8) of \citealp{pmlr-v97-jaber19a}).   
 
 Consider applying this strategy to \mpdag{} $\g$ in Figure \ref{fig:example_meekgraph11}(a), with $\mb{X} = \{X\}$ and $\mb{Y} = \{Y_1,Y_2 \}$. 
Note that $\PossAn(\mb{Y}, \g_{\mb{V} \setminus \mb{X}}) = \{V_1, Y_1, Y_2\}$, that is, $\PossAn(\mb{Y}, \g_{\mb{V} \setminus \mb{X}}) = \mb{V} \setminus \mb{X}$. Then, $\mb{V} \setminus \PossAn(\mb{Y}, \g_{\mb{V} \setminus \mb{X}}) = \mb{X}$. The strategy of \cite{pmlr-v97-jaber19a} would dictate that we can identify the causal effect of $\mb{X}$ on $\mb{Y}$ by first identifying the causal effect of $\mb{X}$ on  $ \mb{V} \setminus \mb{X}$ in $\g$. As we have seen in Example \ref{ex:motivate}, the causal effect of  $\mb{X}$ on  $ \mb{V} \setminus \mb{X}$ in $\g$ is not identifiable, whereas the causal effect of $\mb{X}$ on $\mb{Y}$ is identifiable in $\g$.  Therefore, the approach of \cite{pmlr-v97-jaber19a} is not suitable for general \mpdag{}s. The above counter example arises as a consequence a partially directed cycle in the MPDAG. Hence, a modified approach of \cite{pmlr-v97-jaber19a} may lead to a necessary causal identification algorithm in MPDAGs without partially directed cycles.

A natural question of interest is whether a similar approach to ours can be applied to PAGs. 
Another topic for future work is developing a complete identification formula for conditional causal effects in \mpdag{}s.

\subsubsection*{Acknowledgements}

We  thank S{\o}ren Wengel Mogensen, F.\ Richard Guo, Vincent Roulet, and the reviewers for their insightful comments and suggestions.

\appendix

\section{PRELIMINARIES}

\textbf{Subsequences And Subpaths.} A \textit{subsequence} of a path $p$ is obtained by deleting some nodes from $p$ without changing the order of the remaining nodes. For a path $p = \langle X_1,X_2,\dots,X_m \rangle$, the \textit{subpath} from $X_i$ to $X_k$ ($1\le i\le k\le m)$ is the path $p(X_i,X_k) = \langle X_i,X_{i+1},\dots,X_{k}\rangle$.

\textbf{Concatenation.} We denote concatenation of paths by $\oplus$, so that for a path $p = \langle X_1,X_2,\dots,X_m \rangle$, $p = p(X_1, X_r) \oplus p(X_r, X_m)$, for $1\le r\le m$.

\textbf{D-separation.} If $\mathbf{X}$ and  $\mathbf{Y}$ are d-separated given $\mathbf{Z}$ in a $\DAG$ $\g[D]$, we write $\mathbf{X} \perp_{\g[D]} \mathbf{Y} | \mathbf{Z}$.

\textbf{Possible Descendants.} If there is a possibly causal path from $X$ to $Y$, then $Y$ is a \textit{possible descendant} of $X$. We use the convention that every node is a possible descendant of itself.
The set of possible descendants of $X$ in~$\g$ is $\PossDe(X,\g)$. For a set of nodes $\mathbf{X} \subseteq \mathbf{V}$, we let  $\PossDe(\mathbf{X},\g) = \cup_{X \in \mathbf{X}})  \PossDe(X,\g)$.

\textbf{Bayesian And Causal Bayesian Networks.}
If a density $f$ over $\mb{V}$ is  consistent with $\DAG$ $\g[D] =(\mathbf{V},\mathbf{E})$, then $(\g[D], f)$ form a \textit{Bayesian network}. 
Let $\mb{F}$ be a set of density functions made up of all interventional densities $f(\mb{v'}| do(\mb{x}))$ for any $\mb{X} \subset \mb{V}$ and $\mb{V'} = \mb{V} \setminus \mb{X}$ that are consistent with $\g[D]$ ($\mb{F}$ also includes all observational densities consistent with $\g[D]$), then $(\g[D], \mb{F})$ form a \textit{causal Bayesian network}.

\paragraph{Rules Of The Do-calculus \citep{Pearl2009}.}

Let $\mathbf{X,Y,Z}$ and $\mb{W}$ be pairwise disjoint (possibly empty) sets of nodes in a DAG $\g[D] = (\mb{V,E})$
Let $\g[D]_{\overline{\mathbf{X}}}$ denote the graph obtained by deleting all edges into $\mathbf{X}$ from $\g[D]$. Similarly, let $\g[D]_{\underline{\mathbf{X}}}$ denote the graph obtained by deleting all edges out of $\mathbf{X}$ in~$\g[D]$ and let $\g[D]_{\overline{\mathbf{X}}\underline{\mathbf{Z}}}$ denote the graph obtained by deleting all edges into $\mathbf{X}$ and all edges out of $\mathbf{Z}$ in~$\g[D]$. Let $(\g[D],\mb{F})$ be a causal Bayesian network, the following rules hold for densities in $\mb{F}$.

\textbf{Rule 1} (Insertion/deletion of observations). If $\mathbf{Y} \perp_{\g[D]_{\overline{\mathbf{X}}}} \mathbf{Z} | \mathbf{X} \cup \mathbf{W}$, then
\begin{align}
f(\mathbf{y} | do(\mathbf{x}),\mathbf{w}) = f(\mathbf{y} | do(\mathbf{x}),\mathbf{z,w}). \label{rule1}
\end{align}

\textbf{Rule 2.}  If $\mathbf{Y}  \perp_{\g[D]_{\overline{\mathbf{X}}\underline{\mathbf{Z}}}} \mathbf{Z} | \mathbf{X} \cup \mathbf{W}$, then
\begin{align} \label{rule2}
f(\mathbf{y} | do(\mathbf{x}),do(\mathbf{z}),\mathbf{w}) & = f(\mathbf{y} | do(\mathbf{x}),\mathbf{z,w}).
\end{align}

\textbf{Rule 3.} If $\mathbf{Y} \perp_{\g[D]_{\overline{\mathbf{XZ(W)}}}} \mathbf{Z} | \mathbf{X}  \cup \mb{W}$, then
\begin{align} \label{rule3}
\begin{split}
 f(\mathbf{y} | do(\mathbf{x}), \mb{w}) & = f(\mathbf{y} | do(\mathbf{x}),do(\mathbf{z}),\mb{w}),
\end{split}
\end{align}
where $\mb{Z(W)} =  \mb{Z} \setminus \An(\mb{W}, \g[D]_{\overline{\mb{X}}})$.

\subsection{EXISTING RESULTS}

\begin{theorem}[\textbf{Wright's rule} \citealp{wright1921correlation}] Let $\mathbf{X} = \mathbf{A}\mathbf{X} + \mathbf{\epsilon}$, where $\mathbf{A} \in \mathbb{R}^{k \times  k}$, $\mathbf{X}$ $= (X_1,\dots, X_k)^T$  and $\mathbf{\epsilon} = (\epsilon_1,\dots, \epsilon_k)^T$ is a vector of mutually independent errors with means zero. Moreover, let $\Var(\mathbf{X}) = \mathbf{I}$.
Let $\g[D] = (\mathbf{X},\mathbf{E})$, be the corresponding $\DAG$ such that $X_i \rightarrow X_j$ in $\g[D]$ if and only if $A_{ji} \neq 0$. A nonzero entry $A_{ji}$ is called the edge coefficient of $X_i \rightarrow X_j$.
For two distinct nodes $X_i,X_j \in \mathbf{X}$, let $p_1, \dots, p_r$ be all paths between $X_i$ and $X_j$ in $\g[D]$ that do not contain a collider. Then $\Cov(X_i,X_j) = \sum_{s=1}^{r}\pi_s$, where $\pi_s$ is the product of all edge coefficients along path $p_s$, $s \in \{1,\dots, r\}$.
\label{theorem:wright}
\end{theorem}
\begin{theorem}[c.f.\ Theorem 3.2.4 \citealp{mardia1980multivariate}]
Let $\mathbf{X} = (\mathbf{X_1}^T,\mathbf{X_2}^T)^T$ be a $p$-dimensional multivariate Gaussian random vector with mean vector $\mathbf{\mu} = (\mathbf{\mu_1}^T,\mathbf{\mu_2}^T)^T$ and covariance matrix $\mathbf{\Sigma} = \begin{bmatrix}
\mathbf{\Sigma_{11}} & \mathbf{\Sigma}_{12} \\
\mathbf{\Sigma_{21}} & \mathbf{\Sigma}_{22}
\end{bmatrix}$, so that $\mathbf{X_1}$ is a $q$-dimensional multivariate Gaussian random vector with mean vector $\mathbf{\mu_1}$ and covariance matrix $\mathbf{\Sigma}_{11}$ and  $\mathbf{X_2}$ is a $(p-q)$-dimensional multivariate Gaussian random vector with mean vector $\mathbf{\mu_2}$ and covariance matrix $\mathbf{\Sigma}_{22}$.
Then $E[\mathbf{X_2} | \mathbf{X_1 = x_1}] = \mathbf{\mu_2} + \mathbf{\Sigma}_{21}\mathbf{\Sigma}_{11}^{-1} (\mathbf{x_1} - \mathbf{\mu_1})$.
\label{theorem:mardia-condexp}
\end{theorem}

\begin{lemma}
        Let $\mathbf{X}$ and $\mathbf{Y}$ be disjoint node sets in a \mpdag{} $\g$. Suppose that there is a proper possibly causal path from $\mb{X}$ to $\mb{Y}$ that starts with an undirected edge in $\g$, then        
         there is one such path $q = \langle X, V_1, \dots, Y \rangle$, $X \in \mb{X}$, $Y \in \mb{Y}$ in $\g$ and DAGs $\g[D]^{1}, \g[D]^{2}$ in $[\g]$ such that the path in  $\g[D]^{1}$ consisting of the same sequence of nodes  as $q$ is of the form $X \rightarrow V_1 \rightarrow \dots \rightarrow Y$ and in $\g[D]^{2}$ the path consisting of the same sequence of nodes  as $q$ is of the form $X \leftarrow V_1 \rightarrow \dots \rightarrow Y$.
   \label{lemma:adj amen gen}
\end{lemma}

\begin{lemma}[Lemma 3.2  of \citealp{perkovic17}]
Let $\pstar$ be a path from $X$ to $Y$ in a \mpdag{} $\g$. If $\pstar$ is non-causal in $\g$, then for every $\DAG$ $\g[D]$ in $[\g]$  the corresponding path to $\pstar$  in $\g[D]$ is non-causal. Conversely, if $p$ is a causal path in at least one $\DAG$ $\g[D]$ in $[\g]$, then the corresponding path to $p$ in $\g$ is possibly causal.
\label{lemma:poss-dir-path}
\end{lemma}

\begin{lemma}[Lemma 3.5  of \citealp{perkovic17}]
Let $p = \langle V_1, \dots , V_k \rangle$ be a definite status path in a \mpdag{} $\g$. Then $p$ is possibly causal if and only if there is no $V_i \leftarrow V_{i+1}$, for $i \in \{1,\dots, k-1\}$  in $\g$.
\label{lemma:def-stat-poss-dir}
\end{lemma}

\begin{lemma}[Lemma 3.6 of \citealp{perkovic17}]
 Let $X$ and $Y$ be distinct nodes in a \mpdag{} $\g$. If $p$ is a possibly causal path from $X$ to $Y$ in $\g$, then a subsequence $\pstar$ of $p$ forms a possibly causal unshielded path from $X$ to $Y$ in~$\g$.
\label{lemma:unshielded-analog}
\end{lemma}

\setcounter{algocf}{1}
\begin{algorithm}[t]
\vspace{.2cm}
 \Input{DAG or CPDAG $\mathcal{G} =(\mathbf{V,E})$.}
 \Output{An ordered list $\mb{B} = (\mb{B_1}, \dots , \mb{B_k}), k \ge 1$ of the bucket decomposition of $\mb{V}$ in $\g$.}
 \vspace{.2cm}
   \SetAlgoLined
 Let $\mb{ConComp}$ be the bucket decomposition of $\mb{V}$ in $\g$;\\
   Let $\mb{B}$ be an empty list;\\
\While{$\mb{ConComp}\neq \emptyset$}{
Let  $\mb{C} \in \mb{ConComp}$;\\
Let $\overline{\mb{C}}$ be the set of nodes in $\mb{ConComp}$ that are not in $\mb{C}$;\\
\If{all edges between  $\mb{C}$ and  $\overline{\mb{C}}$   are into $\mb{C}$ in $\g$}
{ Add $\mb{C}$ to the beginning of $\mb{B}$;}
}
\Return  $\mb{B}$;\\
\caption{PTO algorithm  \citep{jaber2018graphical}}
\label{alg:pto1}
\end{algorithm}

\begin{lemma}[c.f.\ Lemma 1 of \citealp{jaber2018graphical}] \label{lemma:1}
 Let $\g = (\mb{V,E})$ be a CPDAG or DAG and let $\mb{B}  = (\mb{B_1}, \dots, \mb{B_k})$, $k \ge 1$, be the output of \texttt{PTO}$(\g)$ (Algorithm \ref{alg:pto1}). Then for each $i,j \in \{1, \dots k\}$, $\mb{B_i}$ and $\mb{B_j}$ are buckets in $\mb{V}$ and if $i < j$, then $\mb{B_i} < \mb{B_j}$.
\end{lemma}

\begin{lemma}[c.f.\ Lemma E.6 of \citealp{henckel}]
 Let $\mb{X}$ and $\mb{Y}$ be disjoint node sets in an MPDAG $\g$ and suppose that there is no proper possibly causal path from $\mb{X}$ to $\mb{Y}$ that starts with an undirected edge in $\g$. Let $\g[D]$ be a DAG in $[\g]$. Then $\fb{\g} \subseteq \De(\mb{X},\g)$.
\label{lemma:forbleo}
\end{lemma}
 
\begin{lemma} \citep[Lemma~A.7 in][]{ernestroth2016} 
Let $X$ and $Y$ be nodes in an  MPDAG $\g = (\mathbf{V,E})$ such that $X - Y$ is in $\g$. Let $\g'$ be an MPDAG constructed from $\g$ by adding $X \to Y$ to $\g$ and completing the orientation rules R1 - R4 of \cite{meek1995causal}. For any $Z,W \in \mathbf{V}$ if $Z \rightarrow W$ is in $\g'$ and $Z -W$ is in $\g$, then $W \in \De(Y,\g')$.
\label{lemma:always-desc}
\end{lemma}

\begin{lemma} \citep[cf.\ Lemma~A.8 in][]{ernestroth2016}
Let $X$ be a node in an  MPDAG $\g  = (\mathbf{V,E})$. Then there is a DAG $\g[D]$,  $\g[D] \in [\g]$ such that $X \rightarrow S$ is in $\g[D]$ for all $X - S$ in $\g$. 
\label{lemma:no-new-into}
\end{lemma}

\section{PROOFS FOR SECTION \ref{sec:id-necess} OF THE MAIN TEXT}

\begin{proofof}[Lemma \ref{lemma:adj amen gen}]
The published version of this paper does not include this proof. Instead, Lemma \ref{lemma:adj amen gen} is stated and references to Lemma C.1 of \citealp{perkovic17}, and Lemma 8 of \cite{perkovic16} are given. This unfortunately is not correct as the result of Lemma \ref{lemma:adj amen gen} does not directly follow from the proofs of these prior results. I am grateful to Sara LaPlante for pointing out this out.
The subtle difference between the proof of this result and that of  Lemma C.1 of \citealp{perkovic17} is in the last paragraph below. 

Let $\pstar[q] = \langle X , \dots , Y \rangle, X \in \mb{X}, Y \in \mb{Y}$ be a proper possibly  causal path from $\mb{X}$ to $\mb{Y}$ in $\g$ that starts with an undirected edge. 
Furthermore, let $q = \langle X= V_0,V_1,\dots,V_k=Y\rangle, k \ge 1$ be a shortest subsequence of $\pstar[q]$ such that $q$ is also a proper possibly causal path that starts with a undirected edge in~$\g$. 
  
   Suppose first that $q$ is of definite status in $\g$. Let $\g[D]_1$ be a $\DAG$ in $[\g]$ that contains $X \rightarrow V_1$  and let $\g[D]_2$ be a $\DAG$ in $[\g]$ that has no additional edges into $V_1$ as compared to $\g$ (Lemma~\ref{lemma:no-new-into}). Then $X \to V_1 \to \dots \to Y$ is in $\g[D]_1$  and $X \leftarrow V_1 \rightarrow \dots \rightarrow Y$ is in $\g[D]_2$ and we are done. 
   
Otherwise, $q$ is not of definite status in $\g$.
By choice of $q$, then $q(V_1,Y)$ must be unshielded and hence, of definite status in $\g$.
Since $q$ is not of definite status and $q(V_1,Y)$ is of definite status, it follows that $V_1$ is not of definite status on $q$. By choice of $q$, $X - V_1$ is in $\g$. Additionally, since $V_1$ is not of definite status on $q$ and since $q$ is a possibly causal path in $\g$, $V_1 - V_2$ is in $\g$ and $X$ is adjacent to $V_2$ in $\g$.
Moreover, we must have $X \to V_2$, since $X - V_2$ contradicts the choice of $q$, and $X \leftarrow V_2$ contradicts that $q$ is possibly causal in $\g$. 

Let $\g[D]_1$ be a $\DAG$ in $\g$ that has no additional edges into $V_1$ compared to $\g$ (Lemma~\ref{lemma:no-new-into}). Since $q(V_1, Y)$ is of definite status in $\g$, $X \leftarrow V_1 \rightarrow V_2 \rightarrow \dots \rightarrow Y$ is in $\g[D]_1$.

Let $\g'$ be an MPDAG constructed from $\g$ by adding edge orientation $V_1 \to V_2$ and completing the orientation rules R1 - R4 of \cite{meek1995causal}. Then $[\g'] \subset [\g]$.
Since $q(V_1, Y)$ is a path of definite status, $V_1 \to \dots \to Y$ is in $\g'$. 
 Since $X \to V_2$ is in $\g$, we know that $X \notin \De(V_2, \g')$. Therefore, by Lemma \ref{lemma:always-desc},  $V_1 \to X$ is not in $\g'$.  Hence $q$ is of the form $X - V_1 \to \dots \to Y$, or $X \to V_1 \to \dots \to Y$ in $\g'$. If $X \to V_1$ is in $\g'$, then we can choose as DAG $\g[D]_2$ any  DAG in $[\g']$. Alternatively, if $X - V_1$ is in $\g'$, we can let $\g[D]_2$ be any DAG in $[\g']$ with $X \to V_1$. Then, once again,  the path corresponding to $q$ will be of the form $X \to V_1 \to \dots \to Y$ in $\g[D]_2$.
\end{proofof}

\begin{proofof}[Proposition \ref{lemmaHedgePDAG}]
This proof follows a similar reasoning as the proof of Theorem 2 of \cite{shpitser2006identification} and proof of Theorem 57 of \cite{perkovic16}. 

By Lemma \ref{lemma:adj amen gen}, there is a proper possibly causal path $q =   \langle X , V_1 , \dots , Y\rangle $, $k \ge 1$, $ X \in \mb{X}$, $Y \in \mb{Y}$  in $\g$ and DAGs $\g[D]^1$ and $\g[D]^2$ in $[\g]$ such that  $X \rightarrow V_1 \rightarrow \dots \rightarrow Y$ is in $\g[D]^1$ and $X \leftarrow V_1 \rightarrow \dots \rightarrow Y$ is in $\g[D]^2$ (the special case when $k=1$ is $X \leftarrow Y$).

Consider a multivariate Gaussian density over $\mb{V}$ with mean vector zero, constructed using a linear structural causal model (SCM) with Gaussian noise. 
In particular, each random variable $A \in \mathbf{V}$ is a linear combination of its parents in $\g[D]^1$ and a designated Gaussian noise variable $\epsilon_{A}$ with zero mean and a fixed variance. The Gaussian noise variables $\{\mathbf{\epsilon}_{A}: A \in \mathbf{V}\}$, are mutually independent.

We define the SCM such that all edge coefficients except for the ones on $q_1$ are $0$, and all edge coefficients  on $q_1$ are in $(0,1)$ and small enough so that we can choose the residual variances so that the variance of every random variable in $\mathbf{V}$ is 1. 

The  density $f$ of $\mb{V}$ generated in this way is consistent with  $\g[D]^{1}$ and thus, $f$ is also consistent with $\g$ and $\g[D]^{2}$ \citep{lauritzen1990independence}. Moreover, $f$ is  consistent with $\DAG$ $\g[D]^{11}$ that is obtained from $\g[D]^{1}$ by removing all edges except for the ones on $q_1$. Analogously, $f$ is also consistent with $\DAG$ $\g[D]^{21}$ that is  obtained from $\g[D]^{2}$ by removing all edges except for the ones on $q_2$. Hence, let $f_1(\mb{v}) = f(\mb{v})$ and let $f_2(\mb{v})= f(\mb{v}) $. 

Let $f_{1}(\mb{v'} |do(\mb{x}))$ be an interventional density consistent with $\g[D]^{11}$. Similarly let $f_{2}(\mb{v'} |do(\mb{x}))$ be an interventional density consistent with $\g[D]^{21}$. Then $f_{1}(\mb{v'} |do(\mb{x}))$  and $f_{1}(\mb{v'} |do(\mb{x}))$  are also interventional densities consistent with $\g[D]^{1}$ and $\g[D]^{2}$, respectively. 
Now, $f_{1}(\mb{y} |do(\mb{x}))$ is a marginal interventional density of $\mb{Y}$ that can be calculated from the density  $f_{1}(\mb{v'} |do(\mb{x}))$  and the analagous is true for  $f_{2}(\mb{y} |do(\mb{x}))$  and $f_{2}(\mb{v'} |do(\mb{x}))$. 

In order to show that  $f_1(\mb{y} | do(\mb{x})) \neq f_2 (\mb{y} | do(\mb{x}))$, it suffices to show that  $f_1(y | do(\mb{x} = 1)) \neq  f_2 (y | do(\mb{x} =\mb{1}))$ for at least one $Y \in \mb{Y}$ when all $\mb{X}$ variables are set to $1$ by a do-intervention. In order for   $f_1(y | do(\mb{x} = 1)) \neq  f_2 (y | do(\mb{x} =\mb{1}))$ to hold, it is enough to show that  the expectation of $Y$ is not the same under these two densities. Hence, let $E_{1}[ Y \mid do(\mathbf{X =1}) ]$ denote the expectation of $Y$, under  $f_1(y | do(\mb{X} = \mb{1}))$ and let  $E_{2}[ Y \mid do(\mathbf{X =1}) ]$   denote the expectation of $\mb{Y}$, under  $f_2(y | do(\mb{X} = \mb{1}))$.

Since $Y$ is d-separated from $\mb{X}$ in $\g[D]^{21}_{\overline{\mb{X}}}$ we can use Rule 3 of the do-calculus (see equation \eqref{rule3}) to conclude that  $E_{2}[ Y \mid do(\mathbf{X}=\mathbf{1})]= E[Y]= 0$. 
Similarly, since $Y$ is d-separated from $\mb{X}$ in $\g[D]^{11}_{\underline{\mb{X}}}$, we can use  Rule 2 of the do-calculus (see equation \eqref{rule2}) to conclude that $E_{1}[ Y \mid do(\mathbf{X}=\mathbf{1})]= E[Y|X=1]$. By Theorems \ref{theorem:mardia-condexp} and \ref{theorem:wright},  $ E[ Y \mid X=1]= \Cov(X,Y) = a$, where $a$ is the product of all edge coefficients on $q_1$. Since  $ a \neq 0$, $E_{1}[ Y \mid do(\mathbf{X =1}) ] \neq E_{2}[ Y \mid do(\mathbf{X =1}) ]$.
\end{proofof}

\section{PROOFS FOR SECTION \ref{sec:id-pto} OF THE MAIN TEXT}

\begin{lemma}\label{lemma:algo-completes}
Let $\mb{D}$ be any subset of $\mb{V}$ in \mpdag{} $\g = (\mb{V},\mb{E})$.
Then the call to algorithm $PCO(\mb{D}, \g)$ will complete. Meaning that, at each iteration of the while loop in $PCO(\mb{D}, \g)$ (Algorithm \ref{alg:pto}), there is a bucket $\mb{C}$ among the remaining buckets in $\mb{ConComp}$ (the bucket decomposition of $\mb{V}$) such that all edges between $\mb{C}$ and $\mb{ConComp} \setminus \mb{C}$ are into $\mb{C}$ in $\g$.
\end{lemma}

\begin{proofof}[Lemma \ref{lemma:algo-completes}]
Let $\mb{C_1},\dots, \mb{C_k}$ be the buckets in $\mb{ConComp}$ at some iteration of the while loop in the call to PCO$(\mb{D}, \g)$. Suppose for contradiction that there is no bucket $\mb{C_i}$, $i \in \{1, \dots, k\}$ such that all edges between $\mb{C_i}$ and $\cup_{j=1}^{k} \mb{C_j} \setminus \mb{C_i}$ are into $\mb{C_i}$. We will show that this leads to the conclusion that $\g$ is not acyclic (a contradiction).

Consider a directed graph $\g_{1}$ constructed so that each bucket in $\mb{ConComp}$ represents one node in $\g_{1}$. Meaning, a bucket $\mb{C_i}$, $i\in \{1, \dots , k\}$ is represented by a node $C_i$ in $\g_{1}$.
Also, let $C_i \to  C_j$, $i,j \in \{1, \dots , k\}$, be in $\g_{1}$ if $A \to B$ is in $\g$ and $A \in C_i$, $B \in C_j$.

Since there is no bucket $\mb{C_i}$ in $\mb{ConComp}$ such that all edges between $\mb{C_i}$ and $\cup_{j=1}^{k} \mb{C_j} \setminus \mb{C_i}$ are into $\mb{C_i}$, there is either a directed cycle in $\g_{1}$, or $C_l \to C_r$ and $C_r \to C_l$ is in $\g_{1}$ for some $l,r \in \{1,...,k\}$. For simplicity, we will refer to both previously mentioned cases as directed cycles.

Let us choose one such directed cycle in $\g_1$, that is, let $C_{r_1}\to \dots \to C_{r_m} \to C_{r_1}$, $ 2 \le  m \le  k$, $r_1, \dots , r_m \in \{1, \dots, k\}$, be in $\g_{1}$. Let $A_i \in \mb{C_{r_i}}$ and $B_{i+1} \in \mb{C_{r_{i+1}}}$, for all $i \in \{1, \dots ,m-1\}$, such that $A_i \to B_{i+1}$ is in $\g$. Additionally, let $A_m \in \mb{C_{r_m}}$, and $B_{1} \in \mb{C_{r_1}}$ such that $A_m \to B_1$ is in $\g$.

Since $A_1 \to B_2$ is in $\g$ and $B_2$ and $A_2$ are in the same bucket $\mb{C_{r_2}}$ in $\g$, by  Lemma \ref{lemma:anc-bucket2}, $A_1 \to A_2$. The same reasoning can be applied to conclude that $A_i \to A_{i+1}$ , for all $i \in \{1,...,m-1\}$  and also that $A_m \to A_1$ is in $\g$. Thus, $A_1 \to A_2 \to \dots \to A_m \to A_1$, a directed cycle is in $\g$, a contradiction.
\end{proofof}

\begin{proofof}[Lemma \ref{lemma:PCO}]
 Lemma \ref{lemma:anc-bucket2} and Lemma \ref{lemma:1} together imply that Algorithm \ref{alg:pto1} can be applied to a \mpdag{} $\g$ and also that the output of \texttt{PTO}($\g$) is the same as that of \texttt{PCO}($\mb{V},\g$). Furthermore, \texttt{PTO}($\g$) =   \texttt{PCO}($\mb{V},\g$) = ($\overline{\mb{B_1}}, \dots  ,\overline{\mb{B_r}}$) $r \ge k$, where  for all $i, j\in \{ 1, \dots, r\}$, $\overline{\mb{B_i}}$ and $\overline{\mb{B_j}}$ are buckets in $\mb{V}$ in $\g$,  and if $i<j$, then $\overline{\mb{B_i}} < \overline{\mb{B_j}}$ with respect to $\g$.
 
The statement of the lemma then follows directly from the definition of buckets (Definition \ref{def:max-bucket}) and Corollary \ref{cor:bucket}, since for each $l \in \{1 ,\dots, k\}$, there exists $s \in \{1, \dots , r\}$ such that $\mb{B_l} =  \mb{D} \cap \overline{\mb{B_s}}$ and $(\mb{B_1}, \dots, \mb{B_k})$ is exactly the  output of \texttt{PCO}($\mb{V},\g$). 
\end{proofof}

\begin{lemma}\label{lemma:anc-bucket2}
Let $\mb{B}$ be a bucket in $\mb{V}$ in  \mpdag{} $\g = (\mb{V},\mb{E})$ and let $X \in \mb{V}$, $X \notin \mb{B}$. 
If there is a causal path from $X$ to $\mb{B}$ in $\g$, then for every node $B \in \mb{B}$ there is a causal path from $X$ to $B$ in $\g$.
\end{lemma}

\begin{proofof}[Lemma \ref{lemma:anc-bucket2}]
Let $p$ be a shortest causal path from $X$ to $\mb{B}$ in $\g$. Then $p$ is of the form $X \rightarrow \dots A \rightarrow B$, possibly $X =A$ and $A \notin \mb{B}$.  

Let $B' \in \mb{B}$, $B' \neq B$ and let $q = \langle B = W_1, \dots , W_r = B' \rangle$, $r >1$ be a shortest undirected path from $B$ to $B'$ in $\g$. It is enough to show that there is an edge $A \rightarrow B'$ is in $\g$.

Since $A \rightarrow B - W_2$, by the properties of \mpdag{}s \citep[][see Figure \ref{fig:orientationRules} in the main text]{meek1995causal}, $A \rightarrow W_2$ or $A - W_2$ is in $\g$. Since $A \notin \mb{B}$, $A \rightarrow W_2$ is in $\g$. If $r = 2$, we are done. Otherwise, $A \rightarrow W_2 - W_3 -\dots - W_k$ is in $\g$ and and we can apply the same reasoning as above iteratively  until we obtain  $ A \rightarrow W_k$ is in $\g$.
\end{proofof}

\section{PROOFS FOR SECTION \ref{sec:factorid} OF THE MAIN TEXT}

The proof of Theorem \ref{thm:id-formula} is given in the main text. Here we provide proofs for the supporting results.

\begin{lemma}\label{lemma:big2}
 Let $\mb{X}$ and $\mb{Y}$ be disjoint node sets in $\mb{V}$ in \mpdag{} $\g =(\mb{V,E})$ and  suppose that there is no proper possibly causal path from $\mb{X}$ to $\mb{Y}$ that starts with an undirected edge in $\g$. Further, let $(\mb{B_1}, \dots \mb{B_k}) =$ \texttt{PCO}$(\An(\mb{Y},\g_{\mb{V} \setminus \mb{X}}),\g)$, $k \ge 1$.

\begin{enumerate}[label=(\roman*)]
\item \label{lemma:amen-bucket22} For $i \in \{1, \dots ,k\}$,  there is no proper possibly causal path from $\mb{X}$ to $\mb{B_i}$ that starts with an undirected edge in $\g$.

\item \label{lemma:dorule2ID2} For $i \in \{2, \dots, k\}$, let $\mb{P_i} = (\cup_{j=1}^{i-1}\mb{B_i}) \cap \Pa(\mb{B_i}, \g)$. Then for every DAG $\g[D]$ in $[\g]$ and every interventional density $f$ consistent with $\g[D]$ we have
$$f(\mb{b_i}| \mb{b_{i-1}},\dots , \mb{b_1}, do(\mb{x})) = f(\mb{b_{i}}| \mb{p_i}, do(\mb{x})).$$

\item \label{lemma:dorule3ID2} For $i \in \{2, \dots, k\}$, let $\mb{P_i} = (\cup_{j=1}^{i-1}\mb{B_i}) \cap \Pa(\mb{B_i}, \g)$.
For $i \in \{1, \dots, k\}$, let $\mb{X_{p_i}} = \mb{X} \cap \Pa(\mb{B_i}, \g)$. Then for every DAG $\g[D]$ in $[\g]$ and every interventional density $f$ consistent with $\g[D]$ we have
$$f(\mb{b_{i}}| \mb{p_i}, do(\mb{x})) = f(\mb{b_{i}}| \mb{p_i}, do(\mb{x_{p_i}})).$$
Additionally, $f(\mb{b_1}|do(\mb{x})) =  f(\mb{b_1}|do (\mb{x_{p_1}})).$

\item \label{lemma:dorulepaID2}  For $i \in \{2, \dots, k\}$, let $\mb{P_i} = (\cup_{j=1}^{i-1}\mb{B_i}) \cap \Pa(\mb{B_i}, \g)$.
For $i \in \{1, \dots, k\}$, let $\mb{X_{p_i}} = \mb{X} \cap \Pa(\mb{B_i}, \g)$. Then for every DAG $\g[D]$ in $[\g]$ and every interventional density $f$ consistent with $\g[D]$ we have
$$ f(\mb{b_{i}}| \mb{p_i}, do(\mb{x_{p_i}})) = f(\mb{b_i}|\pa(\mb{b_i},\g)),$$
for values $\pa(\mb{b_i},\g)$ of $\Pa(\mb{b_i},\g)$ that are in agreement with $\mb{x}$.
\end{enumerate}
\end{lemma}

\begin{proofof}[Lemma \ref{lemma:big2}]
\textbf{ \ref{lemma:amen-bucket22}:} Suppose for a contradiction that there is a proper possibly causal path from $\mb{X}$ to $\mb{B_i}$ that starts with an undirected edge in $\g$. 
Let $p =  \langle X,  \dots, B \rangle$, $X \in \mb{X}$, $B \in \mb{B_i}$, be a shortest such path in $\g$.  Then $p$ is unshielded in $\g$ (Lemma \ref{lemma:unshielded-analog}).

Since $B \in \An(\mb{Y}, \g_{\mb{V} \setminus \mb{X}})$ there is a causal path $q$ from $B$ to $\mb{Y}$ in $\g$ that does not contain a node in $\mb{X}$. 
No node other than $B$ is both on $q$ and $p$ (otherwise, by definition $p$ is not possibly causal from $X$ to $B$). Hence, by Lemma \ref{lemma:concat}, $p \oplus q$ is a proper possibly causal path from $\mb{X}$ to $\mb{Y}$ that starts with an undirected edge in $\g$, which is a contradiction.

\textbf{\ref{lemma:dorule2ID2}:}  Let $\mb{N_i} = (\cup_{j=1}^{i-1} \mb{B_{j}}) \setminus \Pa(\mb{B_i}, \g)$.
If $\mb{B_i}~\perp_{\g[D]_{\overline{\mb{X}}}}~\mb{N_i}~|~(\mb{X} \cup \mb{P_i})$, then by Rule 1 of the do calculus: 
 $f(\mb{b_i}| \mb{b_{i-1}},\dots , \mb{b_1}, do(\mb{x})) = f(\mb{b_{i}}| \mb{p_i}, do(\mb{x}))$ (see equation \eqref{rule1}).

Suppose for a contradiction that there is a path from $\mb{B_i}$ to $\mb{N_i}$ that is d-connecting given $\mb{X} \cup \mb{P_i}$ in $\g[D]_{\overline{\mb{X}}}$. Let $p =\langle B_i, \dots , N \rangle$, $B_i \in \mb{B_i}$, $N \in \mb{N_i}$  be a shortest such path.  Let $\pstar$ be the path in $\g$ that consists of the same sequence of nodes as $p$ in $\g[D]_{\overline{\mb{X}}}$.

First suppose that $p$ is of the form $B_i \rightarrow \dots N$. 
Since $B_i \in \mb{B_i}$ and $\mb{N_i} \subseteq (\cup_{j=1}^{i-1} \mb{B_{j}}) $, $p$ is not causal from $B_i$ to $N$ (Lemma \ref{lemma:PCO}).
Hence, let $C$ be the closest collider to $B_i$ on $p$, that is, $p$ has the form $B_i \rightarrow \dots \rightarrow C \leftarrow \dots N$. Since $p$ is d-connecting given  $\mb{X} \cup \mb{P_i}$ in $\g[D]_{\overline{\mb{X}}}$, $C $ must be an ancestor of  $ \mb{P_i}$ in $\g[D]_{\overline{\mb{X}}}$. However, then there is a causal path from $B_i \in \mb{B_i}$ to $\mb{P_i} \subseteq  (\cup_{j=1}^{i-1} \mb{B_{j}})$ which contradicts Lemma \ref{lemma:PCO}.

Next, suppose that $p$ is of the form $B_i \leftarrow A \dots N$, $A \notin \mb{B_i}$.  Since $\Pa(\mb{B_i}, \g) \subseteq (\mb{X} \cup \mb{P_i})$ and since $p$ is d-connecting given $(\mb{X} \cup \mb{P_i})$, $B_i - A$  is in $\g$ and $A \notin (\mb{X} \cup \mb{P_i})$. 

Note that $\pstar$ cannot be undirected, since that would imply that $N \in \mb{B_i}$ and contradict Lemma \ref{lemma:PCO}. Hence,  let $B$ be the closest node to $B_i$ on $\pstar$ such that $\pstar(B, N)$ starts with a directed edge (possibly $B=A$). Then $\pstar$ is either of the form $B_i -A - \dots  - L - B \rightarrow R \dots N $ or of the form $B_i -A - \dots  - L - B \leftarrow R \dots N $. 

Suppose first that $\pstar$ is of the form $B_i -A - \dots  - L - B \rightarrow R \dots N $. Then $B \notin ( \mb{X} \cup \mb{P_i} \cup \mb{B_i})$ otherwise, $p$ is either blocked by $\mb{X} \cup \mb{P_i}$, or a shorter path could have been chosen. 

Let $(\mb{B_1^{'}}, \dots \mb{B_r^{'}}) =$ \texttt{PCO}$(\mb{V}, \g)$, $r \ge k$. Let $l \in \{i , \dots, r\}$ such that  $\mb{B_l^{'}} \cap \mb{B_i} \neq \emptyset$, then $B_i, B \in \mb{B_l^{'}}$ and $N \in (\cup_{j=1}^{l-1} \mb{B_{j}^{'}})$. 
Now consider subpath $p(B,N)$. By Lemma \ref{lemma:PCO}, $p(B, N)$ cannot be causal from $B$ to $N$. Hence, there is a collider on $p(B, N)$ and we can derive the contradiction using the same reasoning as above. 

Suppose next that  $\pstar$ is of the form $B_i -A - \dots  - L - B \leftarrow R \dots N $. Then either $R \rightarrow L$ or $R - L$ is in $\g$ \citep[][see Figure \ref{fig:example_meekgraph11} in the main text]{meek1995causal}. Then $\langle L, R \rangle$ is also an edge in $\g[D]_{\overline{\mb{X}}}$ otherwise, $L$ or $R$ is in $\mb{X}$ and a non-collider on $p$, so $p$ would be blocked by $\mb{X} \cup \mb{P_i}$. 

Hence, $q= p(B_i,L) \oplus \langle L,R \rangle \oplus p(R,N)$ is a shorter path than $p$ in  $\g[D]_{\overline{\mb{X}}}$. If $L$ and $R$ have the same collider/non-collider status on $q$ on $p$, then $q$ is also d-connecting given $\mb{X} \cup \mb{P_i}$, which would contradict our choice of $p$. Hence, the collider/non-collider status of $L$ or $R$, is different on $p$ and $q$. We now discuss  the cases for the change of collider/non-collider status of $L$ and $R$ and derive a contradiction in each. 

Suppose that $L$ is a collider on $q$, and a non-collider on $p$. This implies that $W \rightarrow L \rightarrow B \leftarrow R$ 
is a subpath of $p$ and $L \leftarrow R$ is in $\g[D]_{\overline{\mb{X}}}$. Even though $L$ is not a collider on $p$, $B$ is a 
collider on $p$ and $L \in \An(B, \g[D]_{\overline{\mb{X}}})$. Since $p$ is d-connecting given $\mb{X} \cup \mb{P_i}$, $
\De(B,  \g[D]_{\overline{\mb{X}}}) \cap (\mb{X} \cup \mb{P_i}) \neq \emptyset$. However, then also $\De(L,  
\g[D]_{\overline{\mb{X}}}) \cap (\mb{X} \cup \mb{P_i}) \neq \emptyset$ and $q$ is also d-connecting given $\mb{X} \cup 
\mb{P_i}$ and a shorter path between $\mb{B_i}$ and $\mb{N_i}$ than $p$, which is a contradiction. 

The contradiction can be derived in exactly the same way as above in the case when $R$ is a collider on $q$, and a non-collider on $p$. Since $B \leftarrow R$ is in $\g[D]_{\overline{\mb{X}}}$, $R$ cannot be anything but a non-collider on $q$, so the only case left to consider is if $L$ is a non-collider on $q$ and a collider on $p$. 

For $L$ to be a  non-collider on $q$ and a collider on $p$, $W \rightarrow L \leftarrow B \leftarrow R$ must be a subpath of $p$ and $L \rightarrow R$ should be in $ \g[D]_{\overline{\mb{X}}}$. But then there is a cycle in $ \g[D]_{\overline{\mb{X}}}$, which is a contradiction. 

\textbf{\ref{lemma:dorule3ID2}:}  We will show that  $f(\mb{b_{i}}| \mb{p_i}, do(\mb{x})) = f(\mb{b_{i}}| \mb{p_i}, do(\mb{x_{p_i}}))$. The simpler case, $f(\mb{b_1}|do(\mb{x})) = f(\mb{b_1} |\do(\mb{x_{p_1}})$ follows from the same proof, when $\mb{B_i}$ is replaced by $\mb{B_1}$ and $\mb{P_i}$ is removed. 

Let $\mb{X_{n_i}} = \mb{X} \setminus \Pa(\mb{B_i}, \g) $ and let  $\mb{X_{n_i}^{'}} = \mb{X_{n_i}}  \setminus \An(\mb{P_i}, \g[D]_{\overline{\mb{X_{p_i}}}})$. That is $X \in \mb{X_{n_i}^{'}}$ if  $X \in \mb{X_{n_i}}$ and if there is no causal path from $X$ to $\mb{P_i}$ in $\g[D]$ that does not contain a node in  $\mb{X_{p_i}}$.

Note that $\Pa(\mb{B_i},\g) = \mb{X_{p_i}} \cup \mb{P_i}$.
 By Rule 3 of the do-calculus, for $f(\mb{b_{i}}| \mb{p_i}, do(\mb{x})) = f(\mb{b_{i}}| \mb{p_i}, do(\mb{x_{p_i}}))$ to hold,  it is enough to show that $\mb{B_i}  \perp_{\g[D]_{\overline{\mb{X_{p_i}}\mb{X_{n_i}^{'}}}}} \mb{X_{n_i}} | \Pa(\mb{B_i}, \g)$ (see equation \eqref{rule3}).

Suppose for a contradiction  that there is a d-connecting path from  $\mb{B_i}$ to $\mb{X_{n_i}}$ in $\g[D]_{\overline{\mb{X_{p_i}}\mb{X_{n_i}^{'}}}}$. Let $p = \langle B_i, \dots , X \rangle$, $B_i \in \mb{B_i}$, $X \in \mb{X_{n_i}}$, be a shortest such path in  $\g[D]_{\overline{\mb{X_{p_i}}\mb{X_{n_i}^{'}}}}$.  Let $\pstar$ be the path in $\g$ that consists of the same sequence of nodes as $p$ in $\g[D]_{\overline{\mb{X_{p_i}}\mb{X_{n_i}^{'}}}}$. This proof follows a very similar line of reasoning to the proof of \ref{lemma:dorule2ID2} above.

Let $(\mb{B_1^{'}}, \dots \mb{B_r^{'}}) =$ \texttt{PCO}$(\mb{V}, \g)$, $r \ge k$. Let $l \in \{i , \dots, r\}$ such that  $\mb{B_l^{'}} \cap \mb{B_i} \neq \emptyset$, then $B_i \in \mb{B_l^{'}}$ and   $ \Pa(\mb{B_i}, \g)\subseteq  (\cup_{j=1}^{i-1} \mb{B_{j}})$.

Suppose that $p$ is of the form $B_i \rightarrow \dots X$. If $X \in  \mb{X_{n_i}^{'}} $, then $p$ is not a causal path since $p$ is a path in $\g[D]_{\overline{\mb{X_{p_i}}\mb{X_{n_i}^{'}}}}$. Otherwise, $X \in \An(\mb{P_i}, \g[D]_{\overline{\mb{X_{p_i}}}})$ and so any causal path from $B_i$ to $X$ would need to contain a node in $\mb{X_{p_i}}$ and hence, would be blocked by  $\Pa(\mb{B_i},\g)$. Thus, $p$ is not a causal path from $B_i$ to $X$. 

Hence,  let $C$ be the closest collider to $B_i$ on $p$, that is, $p$ has the form $B_i \rightarrow \dots \rightarrow C \leftarrow \dots X$. Since $p$ is d-connecting given  $\Pa(\mb{B_i},\g)$, $C $ is be an ancestor of  $\Pa(\mb{B_i},\g)$ in $\g[D]_{\overline{\mb{X_{p_i}}\mb{X_{n_i}^{'}}}}$. However, this would imply that there is a causal path from $B_i \in \mb{B_l^{'}}$ to   $ \Pa(\mb{B_i}, \g)\subseteq  (\cup_{j=1}^{i-1} \mb{B_{j}})$ in  $\g[D]_{\underline{\mb{X_{p_i}}}}$, which contradicts Lemma \ref{lemma:PCO}.

Next, suppose that $p$ is of the form $B_i \leftarrow A \dots X$, $A \notin \mb{B_i}$.  Since $p$  is d-connecting given $ \Pa(\mb{B_i}, \g)$, $A \notin  \Pa(\mb{B_i}, \g)$. Hence,  $B_i - A$  is in $\g$. 

Then $A \in \mb{B_l^{'}}$. Note that by \ref{lemma:amen-bucket22} above, $\mb{X} \cap  \mb{B_l^{'}} = \emptyset$, so $\pstar$ is not an undirected path in $\g$.
 Hence,  let $B$ be the closest node to $B_i$ on $\pstar$ such that $\pstar(B, X)$ starts with a directed edge (possibly $B=A$). Then $\pstar$ is either of the form $B_i -A - \dots  - L - B \rightarrow R \dots X $ or of the form $B_i -A - \dots  - L - B \leftarrow R \dots X $. 

Suppose first that $\pstar$ is of $B_i -A - \dots  - L - B \rightarrow R \dots X $.  Then $B \in \mb{B_l^{'}}$ and so  $B \notin  \mb{X}$. Since $p$ is d-connecting given $\Pa(\mb{B_i}, \g)$, $B \notin \Pa(\mb{B_i}, \g)$ and additionally, $B \notin  \mb{B_i}$ otherwise,  a shorter path could have been chosen. 

Now consider subpath $p(B,X)$. There is at least one collider on $p(B,X)$. Since $B,B_i \in \mb{B_l^{'}}$, the same reasoning as above can be used to derive a contradiction in this case.

Suppose next that  $\pstar$ is of the form $B_i -A - \dots  - L - B \leftarrow R \dots X $. Then either $R \rightarrow L$ or $R - L$ is in $\g$ \citep[][see Figure \ref{fig:example_meekgraph11} in the main text]{meek1995causal}. 
We first show that in either case, edge $\langle L, R\rangle$ is also in $\g[D]_{\overline{\mb{X_{p_i}}\mb{X_{n_i}^{'}}}}$.

Since $L \in \mb{B_l^{'}}$ and since $\mb{X} \cap  \mb{B_l^{'}} = \emptyset$, $L \notin \mb{X}$. Hence, if $R \rightarrow L$ is in $\g$, $R \rightarrow L$ is in  $\g[D]_{\overline{\mb{X_{p_i}}\mb{X_{n_i}^{'}}}}$. If $R - L$ is in $\g$, then $R \in \mb{B_l^{'}}$ and since $\mb{X} \cap  \mb{B_l^{'}} = \emptyset$, $R \notin \mb{X}$, so  $\langle L, R\rangle$ is in $\g[D]_{\overline{\mb{X_{p_i}}\mb{X_{n_i}^{'}}}}$. 

Hence, $q= p(B_i,L) \oplus \langle L,R \rangle \oplus p(R,X)$ is a shorter path than $p$ in $\g[D]_{\overline{\mb{X_{p_i}}\mb{X_{n_i}^{'}}}}$. If $L$ and $R$ have the same collider/non-collider status on $q$ on $p$, then $q$ is also d-connecting given $\Pa(\mb{B_i}, \g)$, which would contradict our choice of $p$. Hence, the collider/non-collider status of $L$ or $R$, is different on $p$ and $q$. We now discuss  the cases for the change of collider/non-collider status of $L$ and $R$ and derive a contradiction in each. 

Suppose that $L$ is a collider on $q$, and a non-collider on $p$. This implies that $W \rightarrow L \rightarrow B \leftarrow R$ 
is a subpath of $p$ and $L \leftarrow R$ are in $\g[D]_{\overline{\mb{X_{p_i}}\mb{X_{n_i}^{'}}}}$. Even though $L$ is not a collider on $p$, $B$ is a  collider on $p$ and $L \in \An(B, \g[D]_{\overline{\mb{X_{p_i}}\mb{X_{n_i}^{'}}}})$. Since $p$ is d-connecting given $\Pa(\mb{B_i}, \g)$, $
\De(B, \g[D]_{\overline{\mb{X_{p_i}}\mb{X_{n_i}^{'}}}}) \cap \Pa(\mb{B_i}, \g) \neq \emptyset$. However, then also $\De(L,  
\g[D]_{\overline{\mb{X_{p_i}}\mb{X_{n_i}^{'}}}}) \cap \Pa(\mb{B_i}, \g) \neq \emptyset$ and $q$ is also d-connecting given $
\Pa(\mb{B_i}, \g)$ and a shorter path between $\mb{B_i}$ and $\mb{X_{n_i}}$ than $p$, which is a contradiction. 

The contradiction can be derived in exactly the same way as above in the case when $R$ is a collider on $q$, and a non-collider on $p$. Since $B \leftarrow R$ is in $\g[D]_{\overline{\mb{X_{p_i}}\mb{X_{n_i}^{'}}}}$, $R$ cannot be anything but a non-collider on $q$, so the only case left to consider is if $L$ is a non-collider on $q$ and a collider on $p$. 

For $L$ to be a  non-collider on $q$ and a collider on $p$, $W \rightarrow L \leftarrow B \leftarrow R$ must be a subpath of $p$ and $L \rightarrow R$ should be in $\g[D]_{\overline{\mb{X_{p_i}}\mb{X_{n_i}^{'}}}}$. But then there is a cycle in $\g[D]_{\overline{\mb{X_{p_i}}\mb{X_{n_i}^{'}}}}$, which is a contradiction. 

\textbf{ \ref{lemma:dorulepaID2}:}. 
If $\mb{B_i} \perp_{\g[D]_{\underline{\mb{X_{p_i}}}}} \mb{X_{p_i}} | \mb{P_i}$, then $ f(\mb{b_{i}}| \mb{p_i}, do(\mb{x_{p_i}})) = f(\mb{b_i}|pa(\mb{b_i},\g))$ by Rule 2 of the do calculus (equation \eqref{rule2}).

Suppose for a contradiction  that there is a d-connecting path from  $\mb{B_i}$ to $\mb{X_{p_i}}$ in $\g[D]_{\underline{\mb{X_{p_i}}}}$. Let $p = \langle B_i, \dots , X \rangle$, $B_i \in \mb{B_i}$, $X \in \mb{X_{p_i}}$, be a shortest such path in  $\g[D]_{\underline{\mb{X_{p_i}}}}$.  Let $\pstar$ be the path in $\g$ that consists of the same sequence of nodes as $p$ in $\g[D]_{\overline{\mb{X}}}$. This proof follows  a very similar line of reasoning to the proof of \ref{lemma:dorule2ID2} above.

Let $(\mb{B_1^{'}}, \dots \mb{B_r^{'}}) =$ \texttt{PCO}$(\mb{V}, \g)$, $r \ge k$. Let $l \in \{i , \dots, r\}$ such that  $\mb{B_l^{'}} \cap \mb{B_i} \neq \emptyset$, then $B_i \in \mb{B_l^{'}}$ and by \ref{lemma:amen-bucket22} above, $\mb{X_{\mb{p_i}}} \subseteq (\cup_{j=1}^{l-1} \mb{B_{j}^{'}})$.

Suppose that $p$ is of the form $B_i \rightarrow \dots X$. Since $B_i \in \mb{B_l^{'}}$ and  $\mb{X_{\mb{p_i}}} \subseteq (\cup_{j=1}^{l-1} \mb{B_{j}^{'}})$, by  Lemma \ref{lemma:PCO}, there is at least one collider on $p$.  Hence, let $C$ be the closest collider to $B_i$ on $p$, that is, $p$ has the form $B_i \rightarrow \dots \rightarrow C \leftarrow \dots X$. Since $p$ is d-connecting given  $ \mb{P_i}$ in $\g[D]_{\underline{\mb{X_{p_i}}}}$, $C $ is be an ancestor of  $ \mb{P_i}$ in $\g[D]_{\underline{\mb{X_{p_i}}}}$. However, this would imply that there is a causal path from $B_i \in \mb{B_i}$ to $\mb{P_i} \subseteq  (\cup_{j=1}^{i-1} \mb{B_{j}})$ in  $\g[D]_{\underline{\mb{X_{p_i}}}}$, which contradicts Lemma \ref{lemma:PCO}.

Next, suppose that $p$ is of the form $B_i \leftarrow A \dots X$, $A \notin \mb{B_i}$.  Since $p$ is a path in $\g[D]_{\underline{\mb{X_{p_i}}}}$, $A \notin \mb{X_{p_i}}$. Additionally, since $p$ is d-connecting given $\mb{P_i}$, $A \notin \mb{P_i}$. Hence,  $B_i - A$  is in $\g$. 

Then $A \in \mb{B_l^{'}}$ and since $X \in (\cup_{j=1}^{l-1} \mb{B_{j}^{'}})$,  $\pstar(A,X)$ is not an undirected path in $\g$.
 Hence,  let $B$ be the closest node to $B_i$ on $\pstar$ such that $\pstar(B, X)$ starts with a directed edge (possibly $B=A$). Then $\pstar$ is either of the form $B_i -A - \dots  - L - B \rightarrow R \dots X $ or of the form $B_i -A - \dots  - L - B \leftarrow R \dots X $. 

Suppose first that $\pstar$ is of $B_i -A - \dots  - L - B \rightarrow R \dots X $.  Then $B \in \mb{B_l^{'}}$ and since $\mb{X_{{p_i}}} \subseteq (\cup_{j=1}^{l-1} \mb{B_{j}^{'}})$,  $B \notin  \mb{X_{p_i}}$. Since $p$ is d-connecting given $\mb{P_i}$, $B \notin \mb{P_i}$ and additionally, $B \notin  \mb{B_i}$ otherwise,  a shorter path could have been chosen. 

Now consider subpath $p(B,X)$. Since $B,B_i \in \mb{B_l^{'}}$, the same reasoning as above can be used to derive a contradiction in this case.

Suppose next that  $\pstar$ is of the form $B_i -A - \dots  - L - B \leftarrow R \dots X $. Then either $R \rightarrow L$ or $R - L$ is in $\g$ \citep[][see Figure \ref{fig:example_meekgraph11} in the main text]{meek1995causal}. Since $R \rightarrow B$ is in $\g[D]_{\underline{\mb{X_{p_i}}}}$, $R \notin \mb{X_{p_i}}$. Since $L \in \mb{B_l^{'}}$, $L \notin \mb{X_{p_i}}$, so  $\langle L, R \rangle$ is also in $\g[D]_{\underline{\mb{X_{p_i}}}}$.

Hence, $q= p(B_i,L) \oplus \langle L,R \rangle \oplus p(R,X)$ is a shorter path than $p$ in $\g[D]_{\underline{\mb{X_{p_i}}}}$. If $L$ and $R$ have the same collider/non-collider status on $q$ on $p$, then $q$ is also d-connecting given $ \mb{P_i}$, which would contradict our choice of $p$. Hence, the collider/non-collider status of $L$ or $R$, is different on $p$ and $q$. We now discuss  the cases for the change of collider/non-collider status of $L$ and $R$ and derive a contradiction in each. 

Suppose that $L$ is a collider on $q$, and a non-collider on $p$. This implies that $W \rightarrow L \rightarrow B \leftarrow R$ 
is a subpath of $p$ and $L \leftarrow R$ are in $\g[D]_{\underline{\mb{X_{p_i}}}}$. Even though, $L$ is not a collider on $p$, $B$ is a  collider on $p$ and $L \in \An(B, \g[D]_{\underline{\mb{X_{p_i}}}})$. Since $p$ is d-connecting given $ \mb{P_i}$, $
\De(B, \g[D]_{\underline{\mb{X_{p_i}}}}) \cap  \mb{P_i} \neq \emptyset$. However, then also $\De(L,  
\g[D]_{\underline{\mb{X_{p_i}}}}) \cap \mb{P_i} \neq \emptyset$ and $q$ is also d-connecting given $
\mb{P_i}$ and a shorter path between $\mb{B_i}$ and $\mb{X_{p_i}}$ than $p$, which is a contradiction. 

The contradiction can be derived in exactly the same way as above in the case when $R$ is a collider on $q$, and a non-collider on $p$. Since $B \leftarrow R$ is in $\g[D]_{\underline{\mb{X_{p_i}}}}$, $R$ cannot be anything but a non-collider on $q$, so the only case left to consider is if $L$ is a non-collider on $q$ and a collider on $p$. 

For $L$ to be a  non-collider on $q$ and a collider on $p$, $W \rightarrow L \leftarrow B \leftarrow R$ must be a subpath of $p$ and $L \rightarrow R$ should be in $\g[D]_{\underline{\mb{X_{p_i}}}}$. But then there is a cycle in $\g[D]_{\underline{\mb{X_{p_i}}}}$, which is a contradiction. 
\end{proofof}

\begin{lemma}\label{lemma:concat}
Let $X,Y$ and $Z$ be distinct nodes in \mpdag{} $\g = (\mb{V,E})$.
Suppose that there is an unshielded possibly causal path $p$ from $X$ to $Y$  and a causal path $q$ from $Y$ to $Z$ in $\g$ such that the only node that $p$ and $q$ have in common is $Y$. 
Then $p \oplus q$ is a possibly causal path from $X$ to $Z$. 
\end{lemma}

\begin{proofof}[Lemma \ref{lemma:concat}]
Suppose for a contradiction that there is an edge $V_q \rightarrow V_p$, where $V_q$ is a node on $q$ and $V_p$ is a node on $p$ (additionally, $V_p \neq Y \neq V_q$). Then $p(V_p,Y)$ cannot be a causal path from $V_p$ to $Y$ since otherwise there is a cycle in $\g$. 
So $p(V_p,Y)$ takes the form $V_p - V_{p+1} \dots Y$.

Let $\g[D]$ be a DAG in $[\g]$, that contains $V_p \rightarrow V_{p+1}$.   Since $p(V_p, Y)$ is an unshielded possibly causal path in $\g$, it corresponds to $V_p \rightarrow \dots \rightarrow Y$ in $\g[D]$. Then $V_{q} \rightarrow V_p \rightarrow \dots \rightarrow Y$ and $q(Y,V_q)$ form a cycle in $\g[D]$, a contradiction.
\end{proofof}

\begin{proofof}[Corollary \ref{cor:truncfac}]
The first statement in Corollary \ref{cor:truncfac} follows from the proof of Theorem \ref{thm:id-formula} when replacing $\mb{Y}$ with $\mb{V}$ and $\mb{X}$ with empty set.

For the second statement in Corollary \ref{cor:truncfac}, note that since there are no undirected edges $X-V$ in $\g$, where $X \in \mb{X}$ and $V \in \mb{V'}$, some of the buckets $\mb{V_i}$, $ i \in \{1, \dots, k\}$ in the bucket decomposition of $\mb{V}$ will contain only nodes in $\mb{X}$. Hence, obtaining the bucket decomposition of $\mb{V'} = \mb{V} \setminus \mb{X}$ is  the same as leaving out  buckets $\mb{V_i}$ that contain only nodes in $\mb{X}$ from  $\mb{V_1}, \dots, \mb{V_k}$. The statement then follows from Theorem \ref{thm:id-formula} when taking $\mb{Y} =  \mb{V'}$. 
\end{proofof}

\section{PROOFS FOR SECTION   \ref{sec:adjustment} OF THE MAIN TEXT}

\begin{proofof}[Proposition \ref{prop:adjustment}]
If the causal effect of $X$ on $Y$ is not identifiable in $\g$, by Theorem \ref{thm:id-formula}, there is a proper possibly causal path from $X$ to $Y$ that starts with an undirected edge in $\g$. Then by Theorem \ref{defadjustmentmpdag}, there is no adjustment set relative to ($X,Y$)  in $\g$.

Hence, suppose that  there is no proper possibly causal path from $X$ to $Y$ that starts with an undirected edge in $\g$ and consider $\Pa(X,\g)$. By Theorem \ref{defadjustmentmpdag}, it is enough to show that $\Pa(X,\g)$ satisfies the generalized adjustment criterion relative to $(\mb{X,Y})$.

 If $\g$ is a DAG, $\Pa(X,\g)$ is an adjustment set relative to $(X,Y)$ by Theorem 3.3.2 of \cite{Pearl2009}.
Hence, suppose that $\g$ is not a DAG. 

Since $\g$ is acyclic, $\Pa(X,\g) \cap \De(X,\g)= \emptyset$. Additionally, by Lemma \ref{lemma:forbleo}, $\f{\g} \subseteq \De(X,\g)$. Hence, $\Pa(X,\g)$ satisfies $\Pa(X,\g) \cap \f{\g} = \emptyset$, that is, condition \ref{cond1} in Theorem \ref{defadjustmentmpdag} relative to ($X,Y$) in $\g$. 

Consider a non-causal definite status path $p$ from $X$ to $Y$. If $p$ is of the form $X \leftarrow \dots Y$ in $\g$, then $p$ is blocked by $\Pa(X,\g)$. If $p$ is of the form $X \rightarrow \dots Y$, then $p$ contains at least one collider  $C \in \De(X,\g)$ and since $\Pa(X,\g) \cap \De(X,\g) = \emptyset$, $p$ is blocked by $\Pa(X,\g)$. 

Lastly, suppose that $p$ is of the form $X - \dots Y$.  Since $p$ is a non-causal path from $X$ to $Y$ and since $p$ is of definite status in $\g$, by Lemma \ref{lemma:def-stat-poss-dir}, there is at least one edge pointing towards $X$ on $p$. Let $D$ be the closest node to $X$  on $p$ such that $p(D,Y)$ is of the form $D \leftarrow \dots Y$ in $\g$. Then by Lemma \ref{lemma:def-stat-poss-dir}, $p(X,D)$ is a possibly causal path from $X$ to $D$ so let $p'$ be an unshielded subsequence of $p(X,D)$ that forms a possibly causal path from $X$ to $D$ in $\g$ (Lemma \ref{lemma:unshielded-analog}).
Additionally, $p$ is of definite status, so $D$ must be a collider on $p$.

In order for $p$ to be blocked by $\Pa(X,\g)$ it is enough to show that $\De(D,\g) \cap \Pa(X,\g) = \emptyset$. Suppose for a contradiction that $E \in \De(D,\g) \cap \Pa(X,\g)$. Let $q$ be  a directed path from $D$ to $E$ in $\g$. Then $p'$ and $q$ satisfy Lemma \ref{lemma:concat} in $\g$, so $p' \oplus q$ is a possibly causal path from $X$ to $E$. By definition of a possibly causal path in \mpdag{}s, this contradicts that $E \in \Pa(X,\g)$.
\end{proofof}

\begin{lemma}\label{lemma:preproc}
Let $\mb{X}$ and $\mb{Y}$ be disjoint node sets in an  \mpdag{} $\g = (\mb{V},\mb{E})$. 
If   there is no possibly causal path from $\mb{X}$ to $\mb{Y}$ in $\g$, then for any  observational density $f$  consistent with $\g$ we have
\begin{align*} 
f(\mb{y}|do(\mb{x})) =f(\mb{y}).
\end{align*}
\end{lemma}

\begin{proofof}[Lemma \ref{lemma:preproc}]
Lemma \ref{lemma:preproc} follows from Lemma \ref{lemma:poss-dir-path}  and Rule 3 of the do-calculus of \cite{Pearl2009} (see equation \eqref{rule3}).  
\end{proofof}

\vskip 0.2in
\bibliographystyle{apalike}
\bibliography{biblioteka}
\end{document}